\documentclass[12pt]{amsart}
\usepackage{amsmath,amsfonts,amsthm,amscd}
\newtheorem{theorem}{Theorem}
\newtheorem{claim}[theorem]{Claim}
\newtheorem{proposition}[theorem]{Proposition}
\newtheorem{lemma}[theorem]{Lemma}

\newtheorem{corollary}[theorem]{Corollary}
\theoremstyle{remark}
\newtheorem*{remark}{Remark}
\newcommand{\ds}{\displaystyle}

\numberwithin{equation}{section} \numberwithin{theorem}{section}
\begin{document}
\bibliographystyle{amsalpha}
\title[Ricci tensor]{Prescribing symmetric functions of the eigenvalues
of the Ricci tensor}
\author{Matthew J. Gursky}
\address{Matthew J. Gursky\\
 Department of Mathematics \\ University of Notre Dame \\ Notre Dame, IN 46556}
\email{mgursky@nd.edu}
\thanks{The research of the first author was
partially supported by NSF Grant DMS-0200646.}
\author{Jeff A. Viaclovsky}
\address{Jeff A. Viaclovsky, Department of Mathematics, MIT, Cambridge, MA 02139}
\email{jeffv@math.mit.edu}
\thanks{The research of the second author was partially
supported by NSF Grant DMS-0202477.}
\date{September 10, 2004}
\begin{abstract}We study the problem of conformally deforming a
metric to a prescribed symmetric function of the eigenvalues of the
Ricci tensor. We prove an existence theorem for a wide class of
symmetric functions on manifolds with positive Ricci curvature,
provided the conformal class admits an {\em{admissible}} metric.

\end{abstract}
\maketitle \setcounter{tocdepth}{1} 
\section{Introduction}

Let $(M^n,g)$ be a smooth, closed Riemannian manifold of dimension
$n$.  We denote the Riemannian curvature tensor by $Riem$, the Ricci
tensor by $Ric$, and the scalar curvature by $R$.  In addition, the
{\it Weyl-Schouten tensor} is defined by
\begin{align}
\label{WStensor} A = \frac{1}{(n-2)}\Big( Ric - \frac{1}{2(n-1)} R
g\Big).
\end{align}
This tensor arises as the ``remainder" in the standard decomposition
of the curvature tensor
\begin{align}
\label{rot} Riem = W + A \odot g,
\end{align}
where $W$ denotes the Weyl curvature tensor and $\odot$ is the
natural extension of the exterior product to symmetric
$(0,2)$-tensors (usually referred to as the {\it Kulkarni-Nomizu}
product, \cite {Besse}). Since the Weyl tensor is conformally
invariant, an important consequence of the decomposition (\ref{rot})
is that the tranformation of the Riemannian curvature tensor under
conformal deformations of metric is completely determined by the
transformation of the symmetric $(0,2)$-tensor $A$.

In \cite{Jeff1} the second author initiated the study of the fully
nonlinear equations arising from the transformation of $A$ under
conformal deformations.  More precisely, let $g_{u} =e^{-2u}g$
denote a conformal metric, and consider the equation
\begin{align}
\label{sigmak} \sigma_k^{1/k}(g_u^{-1}A_u) = f(x),
\end{align}
where $\sigma_k : \mathbf{R}^n \to \mathbf{R}$ denotes the
elementary symmetric polynomial of degree $k$, $A_u$ denotes the
Weyl-Schouten tensor with respect to the metric $g_{u}$, and
$\sigma_k^{1/k}(g_u^{-1}A_u)$ means $\sigma_k(\cdot)$ applied to the
eigenvalues of the $(1,1)$-tensor $g_u^{-1}A_u$ obtained by
``raising an index" of $A_u$.  Following the conventions of our
previous paper \cite{GVAM}, we interpret $A_u$ as a bilinear form on
the tangent space with inner product $g$ (instead of $g_u$). That
is, once we fix a background metric $g$, $\sigma_k(A_u)$ means
$\sigma_k(\cdot)$ applied to the eigenvalues of the $(1,1)$-tensor
$g^{-1}A_u$. To understand the practical effect of this convention,
recall that $A_u$ is related to $A$ by the formula
\begin{align}
A_u = A + \nabla^2 u + du \otimes du - \frac{1}{2}| \nabla u|^2 g
\label{Achange}
\end{align}
(see \cite{Jeff1}).  Consequently, (\ref{sigmak}) is equivalent to
\begin{align}
\label{hessk} \sigma_k^{1/k}(A + \nabla^2 u + du \otimes du -
\frac{1}{2}| \nabla u|^2 g) = f(x)e^{-2u}.
\end{align}
Note that when $k=1$, then $\sigma_1(g^{-1}A) = trace(A) =
\frac{1}{2(n-1)}R$.  Therefore, (\ref{hessk}) is the problem of
prescribing scalar curvature.

To recall the ellipticity properties of (\ref{hessk}), following
\cite{Garding} and \cite{CNSIII} we let $\Gamma_k^{+} \subset
\mathbf{R}^n$ denote the component of $\{x \in \mathbf{R}^n |
\sigma_k(x) > 0 \}$ containing the positive cone $\{ x \in
\mathbf{R}^n | x_1 > 0,..., x_n > 0 \}$.  A solution $u \in
C^2(M^n)$ of (\ref{hessk}) is elliptic if the eigenvalues of $A_u$
are in $\Gamma_k^{+}$ at each point of $M^n$; we then say that $u$
is {\it admissible} (or $k$-{\it admissible}).   By a result of the
second author, if $u \in C^2(M^n)$ is a solution of (\ref{hessk})
and the eigenvalues of $A = A_g$ are everywhere in $\Gamma_k^{+}$,
then $u$ is admissible (see \cite{Jeff1}, Proposition 2). Therefore,
we say that a metric $g$ is {\it k-admissible} if the eigenvalues of
$A = A_g$ are in $\Gamma_k^{+}$, and write $g \in
\Gamma_k^{+}(M^n)$.

In this paper we are interested in the case $k > n/2$. According to
a result of Guan-Viaclovsky-Wang \cite{GVW}, a $k$-admissible metric
with $k > n/2$ has positive Ricci curvature; this is the geometric
significance of our assumption.  Analytically, when $k
> n/2$ we can establish an integral estimate for solutions of
(\ref{hessk}) (see Theorem \ref{intgrad}). As we shall see, this
estimate is used at just about every stage of our analysis. Our main
result is a general existence theory for solutions of (\ref{hessk}):

\begin{theorem}  \label{Main}
Let $(M^n,g)$ be closed $n$-dimensional Riemannian manifold, and
assume \vskip.1in \noindent $(i)$ $g$ is $k$-admissible with $k >
n/2$, and \vskip.1in \noindent $(ii)$ $(M^n,g)$ is not conformally
equivalent to the round $n$-dimensional sphere. \vskip.1in  Then
given any smooth positive function $f \in C^{\infty}(M^n)$ there
exists a solution $u \in C^{\infty}(M^n)$ of (\ref{hessk}), and the
set of all such solutions is compact in the $C^m$-topology for any
$m \geq 0$.
\end{theorem}
\begin{remark}The second assumption above is of course necessary, since the set
of solutions of (\ref{hessk}) on the round sphere with $f(x) =
constant$ is non-compact, while for variable $f$ there are
obstructions to existence.  In particular, there is a ``Pohozaev
identity'' for solutions of (\ref{hessk}) which holds in the
conformally flat case; see \cite{Jeff3}. This identity yields
non-trivial Kazdan-Warner-type obstructions to existence (see
\cite{KazdanWarner}) in the case $(M^n,g)$ is conformally equivalent
to $(S^n, g_{round})$. It is an interesting problem to characterize
the functions $f(x)$ which may arise as $\sigma_k$-curvature
functions in the conformal class of the round sphere, but we do not
address this problem here.
\end{remark}

\subsection{Prior results}
  Due to the amount of activity it has become increasingly difficult
to provide even a partial overview of results in the literature
pertaining to (\ref{hessk}). Therefore, we will limit ourselves to
those which are the most relevant to our work here.

In \cite{Jeff2}, the second author established global {\it a priori}
$C^1$- and $C^2$-estimates for $k$-admissible solutions of
(\ref{hessk}) that depend on $C^0$-estimates.  Since (\ref{hessk})
is a convex function of the eigenvalues of $A_u$, the work of Evans
and Krylov (\cite{Evans}, \cite{Krylov}) give $C^{2,\alpha}$ bounds
once $C^2$-bounds are known. Consequently, one can derive estimates
of all orders from classical elliptic regularity, provided $C^0$-
bounds are known.  Subsequently, Guan and Wang (\cite{GuanWang1})
proved local versions of these estimates which only depend on a
lower bound for solutions on a ball.  Their estimates have the added
advantage of being scale-invariant, which is crucial in our
analysis.  For this reason, in Section 2 of the present paper we
state the main estimate of Guan-Wang and prove some straightforward
but very useful corollaries.

Given $(M^n,g)$ with $g \in \Gamma_k^{+}(M^n)$, finding a solution
of (\ref{hessk}) with $f(x) = constant$ is known as the {\it
$\sigma_k$-Yamabe problem}. In \cite{GVAM} we described the
connection between solving the $\sigma_k$-Yamabe problem when $k
> n/2$ and a new conformal invariant called
the {\it maximal volume} (see the Introduction of \cite{GVAM}). On
the basis of some delicate global volume comparison arguments, we
were able to give sharp estimates for this invariant in dimensions
three and four.  Then, using the local estimates of Guan-Wang and
the Liouville-type theorems of Li-Li \cite{LiLi2}, we proved the
existence and compactness of solutions of the $\sigma_k$-Yamabe
problem for any $k$-admissible four-manifold $(M^4,g)$ ($k \geq 2$),
and any simply connected $k$-admissible three-manifold $(M^3,g)$ ($k
\geq 2$).  More generally, we proved the existence of a number
$C(k,n) \geq 1$, such that if the fundamental group of $M^n$
satisfies $\| \pi_1(M^n) \| > C(k,n)$ then the conformal class of
any $k$-admissible admissible metric with $k > n/2$ admits a
solution of the $\sigma_k$-Yamabe problem. Moreover, the set of all
such solutions is compact.

We note that the proof of Theorem \ref{Main} does not rely on the
Liouville theorem of Li-Li.  Indeed, other than the local estimates
of Guan-Wang, the present paper is fairly self-contained.

There are several existence results for (\ref{hessk}) when $(M^n,g)$
is assumed to be locally conformally flat and $k$-admissible.  In
\cite{LiLi2}, Li and Li solved the $\sigma_k$-Yamabe problem for any
$k \geq 1$, and established compactness of the solution space
assuming the manifold is not conformally equivalent to the sphere.
Guan and Wang (\cite{GuanWang2}) used a parabolic version of
(\ref{hessk}) to prove global existence (in time) of solutions and
convergence to a solution of the $\sigma_k$-Yamabe problem. However,
as we observed above, if $(M^n,g)$ is $k-$admissible with $k > n/2$
then $g$ has positive Ricci curvature; by Myer's Theorem the
universal cover $X^n$ of $M^n$ must be compact, and Kuiper's Theorem
implies $X^n$ is conformally equivalent to the round sphere.  We
conclude the manifold $(M^n,g)$ must be conformal to a spherical
space form. Consequently, there is no significant overlap between
our existence result and those of Li-Li or Guan-Wang.

For global estimates the aforementioned result of Viaclovsky
(\cite{Jeff2}) is optimal: since (\ref{sigmak}) is invariant under
the action of the conformal group, {\it a priori} $C^0$-bounds may
fail for the usual reason (i.e., the conformal group of the round
sphere). Some results have managed to distinguish the case of the
sphere, thereby giving bounds when the manifold is not conformally
equivalent to $S^n$. For example, \cite{CGY2} proved the existence
of solutions to (\ref{hessk}) when $k=2$ and $g$ is $2$-admissible,
for any function $f(x)$, provided $(M^4,g)$ is not conformally
equivalent to the sphere. In \cite{Jeff2} the second author studied
the case $k=n$, and defined another conformal invariant associated
to admissible metrics. When this invariant is below a certain value,
one can establish $C^0$-estimates, giving existence and compactness
for the determinant case on a large class of conformal manifolds.

\subsection{Outline of proof}

In this paper our strategy is quite different from the results just
described.  We begin by defining a $1$-parameter family of equations
that amounts to a deformation of (\ref{hessk}).  When the parameter
$t = 1$, the resulting equation is exactly (\ref{hessk}), while for
$t = 0$ the 'initial' equation is much easier to analyze. This
artifice appears in our previous paper \cite{GVAM}, except that here
we are attempting to solve (\ref{hessk}) for general $f$ and not
just $f(x) = constant$.  In both instances the key observation is
that the Leray-Schauder degree, as defined in the paper of Li
\cite{Yanyan2}, is non-zero. By homotopy-invariance of the degree
the question of existence reduces to establishing {\it a priori}
bounds for solutions for $t \in [0,1]$.

To prove such bounds we argue by contradiction. That is, we assume
the existence of a sequence of solutions $\{ u_i \}$ for which a
$C^0$-bound fails, and undertake a careful study of the blow-up. On
this level our analysis parallels the blow-up theory for solutions
of the Yamabe problem as described, for example, in \cite{Schoen1}.

The first step is to prove a kind of weak compactness result for a
sequence of solutions $\{ u_i \}$, which says that there is a finite
set of points $\Sigma = \{x_1,\dots,x_{\ell} \} \subset M^n$ with
the property that the $u_i$'s are bounded from below and the
derivatives up to order two are uniformly bounded on compact subsets
of $M^n \setminus \Sigma$ (see Proposition \ref{nearsigv1}). This
leads to two possibilities: either a subsequence of $\{ u_i \}$
converges to a limiting solution on $M^n \setminus \Sigma$, or $u_i
\to +\infty$ on $M^n \setminus \Sigma$. Using our integral gradient
estimate, we are able to rule out the former possibility.

The next step is to normalize the sequence $\{ u_i \}$ by choosing a
``regular" point $x_0 \notin \Sigma$ and defining $w_i(x) = u_i(x) -
u_i(x_0)$.  By our preceding observations, a subsequence of $\{ w_i
\}$ converges on compact subsets of $M^n \setminus \Sigma$ in
$C^{1,\alpha}$ to a limit $w \in C^{1,1}_{loc}(M^n \setminus
\Sigma)$.  At this point, the analysis becomes technically quite
different from that of the Yamabe problem, where a divergent
sequence (after normalizing in a similar way) is known to converge
off the singular set to a solution of $L\Gamma = 0$, where $\Gamma$
is a linear combination of fundamental solutions of the conformal
laplacian $L = \Delta  - \frac{(n-2)}{4(n-1)}R$.  By contrast, in
our case the limit is only a {\it viscosity} solution of
\begin{align} \label{viscI}
\sigma_k^{1/k}(A + \nabla^2 w + dw \otimes dw - \frac{1}{2}| \nabla
w|^2 g) \geq 0
\end{align}
In addition, we have no {\it a priori} knowledge of the behavior of
singular solutions of (\ref{viscI}). For example, it is unclear what
is meant by a fundamental solution in this context.

Keeping in mind the goal, if not the means of \cite{Schoen1}, we
remind the reader that Schoen applied the Pohozaev identity to the
singular limit $\Gamma$ to show that the constant term in the
asymptotic expansion of the Green's function has a sign, thus
reducing the problem to the resolution of the Positive Mass Theorem.
In other words, analysis of the sequence is reduced to analysis of
the asymptotically flat metric $\Gamma^{4/(n-2)}g$. For example, if
$(M^n,g)$ is the round sphere then the singular metric defined by
the Green's function $\Gamma_{p}$ with pole at $p$ is flat; in fact,
$(M^n \setminus \{p\},\Gamma_p^{4/(n-2)}g)$ is isometric to
$(\mathbf{R}^n, g_{Euc})$.

Our approach is to also study the manifold $(M^n \setminus \Sigma
,e^{-2w}g)$ defined by the singular limit.  However, the metric $g_w
= e^{-2w}g$ is only $C^{1,1}$, and owing to our lack of knowledge
about the behavior of $w$ near the singular set $\Sigma$, initially
we do not know if $g_w$ is complete. Therefore in Section 6 we
analyze the behavior of $w$, once again relying on the integral
gradient estimate and a kind of weak maximum principle for singular
solutions of (\ref{viscI}). Eventually we are able to show that near
any point $x_k \in \Sigma$,
\begin{align} \label{wfromaboveI}
2 \log d_g(x,x_k) - C \leq w(x) \leq 2 \log d_g(x,x_k) + C
\end{align}
for some constant $C$, where $d_g$ is the distance function with
respect to $g$.  If $\Gamma$ denotes the Green's function for $L$
with singularity at $x_k$, then (\ref{wfromaboveI}) is equivalent to
\begin{align*}
c^{-1} \Gamma(x)^{4/(n-2)} \leq e^{-2w(x)} \leq c
\Gamma(x)^{4/(n-2)}
\end{align*}
for some constant $c > 1$.  Thus, the asymptotic behavior of the
metric $g_w$ at infinity is the same--at least to first order--as
the behavior of $\Gamma^{4/(n-2)}g$.  Consequently, $g_w$ is
complete (see Proposition \ref{comp1}).

The estimate (\ref{wfromaboveI}) can be slightly refined; if
$\Psi(x) = w(x) - 2 \log d_g(x,x_k)$ then (\ref{wfromaboveI}) says
$\Psi(x) = O(1)$ near $x_k$.  In fact, we can show that
\begin{align}  \label{lilbetter}
\int_{U_k} |\nabla \Psi|_{g_w}^n dvol_{g_w}  < \infty
\end{align}
for some neighborhood $U_k$ of $x_k$ (see Theorem \ref{gdecay}).
Using this bound we proceed to analyze the manifold $(M^n,g_w)$ near
infinity. First, we observe that since $g_w$ is the limit of smooth
metrics with positive Ricci curvature, by Bishop's theorem the
volume growth of large balls is sub-Euclidean:
\begin{align} \label{Volbd}
\displaystyle \frac{Vol_{g_w}(B(p_0,r))}{r^n} \leq \omega_n,
\end{align}
where $p_0 \in M \setminus \Sigma$ is a basepoint. Moreover, the
ratio in (\ref{Volbd}) is non-increasing as a function of $r$. Also,
using (\ref{lilbetter}) and a tangent cone analysis, we find that
\begin{align*}
\lim_{r \to \infty} \displaystyle \frac{Vol_{g_w}(B(p_0,r))}{r^n} =
\omega_n.
\end{align*}
Therefore, equality holds in (\ref{Volbd}), which by Bishop's
theorem implies that $g_w$ is isometric to the Euclidean metric. We
emphasize that since the limiting metric $g_w$ is only $C^{1,1}$, we
cannot directly apply the standard version of Bishop's theorem; this
problem makes our arguments technically more difficult.  However,
once equality holds in (\ref{Volbd}), it follows that $w$ is
regular, $e^{-2w} = \Gamma^{4/(n-2)}$, and $(M^n,g)$ is conformal to
the round sphere.

Because much of the technical work of this paper is reduced to
understanding singular solutions which arise as limits of sequences,
we are optimistic that our techniques can be used to study more
general singular solutions of (\ref{hessk}), as in the recent work
of Maria del Mar Gonzalez \cite{Gonzalez2},\cite{Gonzalez1}. Also,
the importance of the integral estimate Theorem \ref{intgrad}
indicates that it should be of independent interest in the study of
other conformally invariant fully nonlinear equations.

\subsection{Other symmetric functions}

Our method of analyzing the blow-up for sequences of solutions to
(\ref{hessk}) can be applied to more general examples of symmetric
functions, provided the Ricci curvature is strictly positive and the
appropriate local estimates are satisfied. To make this precise, let
\begin{align} \label{Fex}
F : \Gamma \subset \mathbf{R}^n \rightarrow \mathbf{R}
\end{align}
with $F \in C^{\infty}(\Gamma) \cap C^{0}(\overline{\Gamma})$, where
$\Gamma \subset \mathbf{R}^n$ is an open, symmetric, convex cone.

We impose the following conditions on the operator $F$:

\vskip.1in $(i)$ $F$ is symmetric, concave, and homogenous of degree
one.

\vskip.1in  $(ii)$ $F > 0$ in $\Gamma$, and $F = 0$ on $\partial
\Gamma$.

\vskip.1in  $(iii)$  $F$ is {\em elliptic}: $F_{\lambda_i}(\lambda)
> 0$ for each $1 \leq i \leq n$, $\lambda \in \Gamma$.

\vskip.1in  $(iv)$ $\Gamma \supset \Gamma_n^{+}$, and there exists a
constant $\delta> 0$ such that any $\lambda =
(\lambda_1,\dots,\lambda_n) \in \Gamma$ satisfies
\begin{align} \label{Ricposcond}
\lambda_i > -\frac{(1- 2\delta)}{(n-2)}\big(\lambda_1 + \dots +
\lambda_n \big) \quad \forall \ 1 \leq i \leq n.
\end{align}

\vskip.1in

To explain the significance of (\ref{Ricposcond}), suppose the
eigenvalues of the Schouten tensor $A_g$ are in $\Gamma$ at each
point of $M^n$.  Then $(M^n,g)$ has positive Ricci curvature: in
fact,
\begin{align} \label{Ricposintro}
Ric_g - 2\delta \sigma_1(A_g)g \geq 0.
\end{align}

For $F$ satisfying $(i)-(iv)$, consider the equation
\begin{align} \label{hessF}
F(A_u) = f(x)e^{-2u},
\end{align}
where we assume $A_u \in \Gamma$ (i.e., $u$ is {\em
$\Gamma$-admissible}). Some examples of interest are

\vskip.25in \noindent {\bf Example 1.}  $F(A_u) =
\sigma_k^{1/k}(A_u)$ with $\Gamma = \Gamma_k^{+}$, $k > n/2$. Thus,
(\ref{hessk}) is an example of (\ref{hessF}).

\vskip.25in \noindent {\bf Example 2.} Let $1 \leq l < k$ and $k
> n/2$, and consider
\begin{align} \label{rats}
F(A_u) = \ds \Big( \frac{\sigma_k(A_u)}{\sigma_l(A_u)}
\Big)^{\frac{1}{k-l}}.
\end{align}
In this case we also take $\Gamma = \Gamma_k^{+}$.

\vskip.25in \noindent {\bf Example 3.} For $\tau \leq 1$ let
\begin{align}
\label{WStensort} A^{\tau} = \frac{1}{(n-2)} \Big( Ric -
\frac{\tau}{2(n-1)} R g \Big)
\end{align}
and consider the equation
\begin{align}
\label{sigmakt} F(A_u) = \sigma_k^{1/k}(A^{\tau}_u) = f(x)e^{-2u}.
\end{align}
By (\ref{Achange}), this is equivalent to the fully nonlinear
equation
\begin{align}
\label{hesskt} \sigma_k^{1/k} \Big(
 A^{\tau} + \nabla^2 u + \frac{1-\tau}{n-2}(\Delta u)g
+ du \otimes du - \frac{2-\tau}{2} |\nabla u|^2 g \Big) =
f(x)e^{-2u}.
\end{align}
In the Appendix we show that the results of \cite{GVW} imply the
existence of $\tau_0 = \tau_0(n,k)
> 0$ and $\delta_0 = \delta(k,n)
> 0$ so that if $1 \geq \tau > \tau_0(n,k)$ and  $A_g^{\tau} \in
\Gamma_k$ with $k>n/2$, then $g$ satisfies (\ref{Ricposcond}) with
$\delta = \delta_0.$

\vskip.2in

For the existence part of our proof we use a degree theory argument
which requires us to introduce a $1$-parameter family of auxiliary
equations. For this reason, we need to consider the following
slightly more general equation:

\begin{align} \label{hessFlong}
F(A_u + G(x) ) = f(x)e^{-2u} + c,
\end{align}
where $G(x)$ is a symmetric $(0,2)$-tensor with eigenvalues in
$\Gamma$, and $c \geq 0$ is a constant. To extend our compactness
theory to equations like (\ref{hessFlong}), we need to verify that
solutions satisfy local estimates like those proved by Guan-Wang.
Such estimates were recently proved by S. Chen \cite{Sophie1}:

\begin{theorem}  {\em (See \cite{Sophie1}, Corollary 1)} \label{SophieEstimate}  Let $F$ satisfy the properties $(i)-(iv)$
above.  If $u \in C^4(B(x_0,\rho))$ is a solution of
(\ref{hessFlong}), then there is a constant
\begin{align*}
C_0 = C_0(n,\rho,\|g\|_{C^4(B(x_0,\rho))}, \|f\|_{C^2(B(x_0,\rho))},
\|G\|_{C^2(B(x_0,\rho))}, c)
\end{align*}
such that
\begin{align} \label{SIC2}
|\nabla^2 u|(x) + |\nabla u|^2(x) \leq C_0\big(1 + e^{-2
\inf_{B(x_0,\rho)}u}\big)
\end{align}
for all $x \in B(x_0,\rho/2)$.
\end{theorem}

An important feature of (\ref{SIC2}) is that both sides of the
inequality have the same homogeneity under the natural dilation
structure of equation (\ref{hessFlong}); see the proof of Lemma
\ref{gradrest} and the Remark following Proposition
\ref{epsilonregr}.

We note that higher order regularity for solutions of
(\ref{hessFlong}) will follow from pointwise bounds on the solution
and its derivatives up to order two, by the aforementioned results
of Evans \cite{Evans} and Krylov \cite{Krylov}.  The point here is
that $C^2$-bounds, along with the properties $(i)-(iv)$, imply that
equation (\ref{hessFlong}) is {\it uniformly} elliptic.  Since this
is not completely obvious we provide a proof in the Appendix.

For the examples enumerated above, local estimates have already
appeared in the literature.  As we already noted, Guan and Wang
established local estimates for solutions of (\ref{hessk}) in
\cite{GuanWang1}.  In subsequent papers
(\cite{GuanWang3},\cite{GuanWangquot}) they proved a similar
estimate for solutions of (\ref{rats}).  In \cite{LiLi2}, Li and Li
proved local estimates for solutions of (\ref{hesskt})  (see also
\cite{GVJDG}).  In both cases, the estimates can be adapted to the
modified equation (\ref{hessFlong}) with very little difficulty. The
work of S. Chen, in addition to giving a unified proof of these
results, also applies to other fully nonlinear equations in
geometry. Applying her result, our method gives
\begin{theorem}  \label{Main2}
Suppose $F: \Gamma \rightarrow \mathbf{R}$ satisfies $(i)-(iv)$. Let
$(M^n,g)$ be closed $n$-dimensional Riemannian manifold, and assume
\vskip.1in \noindent $(i)$ $g$ is $\Gamma$-admissible, and
\vskip.1in \noindent $(ii)$ $(M^n,g)$ is not conformally equivalent
to the round $n$-dimensional sphere. \vskip.1in Then given any
smooth positive function $f \in C^{\infty}(M^n)$ there exists a
solution $u \in C^{\infty}(M^n)$ of
\begin{align*}
F(A_u) = f(x)e^{-2u},
\end{align*}
and the set of all such solutions is compact in the $C^m$-topology
for any $m \geq 0$.
\end{theorem}
Note that, in particular, the symmetric functions arising in
Examples $2$ and $3$ above fall under the umbrella of Theorem
\ref{Main2}.
To simplify the exposition in the paper we give the proof of Theorem
\ref{Main}, while providing some remarks along the way to point out
where modifications are needed for proving Theorem \ref{Main2} (in
fact, there a very few).

\subsection{Acknowledgements}
 The authors would like to thank Alice Chang, Pengfei Guan,
Yanyan Li, and Paul Yang for enlightening discussions on conformal
geometry. The authors would also like to thank Luis Caffarelli and
Yu Yuan for valuable discussions about viscosity solutions. Finally,
we would like to thank Sophie Chen for bringing her work on local
estimates to our attention.

\section{The deformation}
\label{deformation}

Let $f \in C^{\infty}(M^n)$ be a positive function, and for $0 \leq
t \leq 1$ consider the family of equations
\begin{align} \label{start} \begin{split}
\sigma_k^{1/k}& \Big( \lambda_k(1 - \psi(t))g + \psi(t)A +
\nabla^2 u + du \otimes du - \frac{1}{2}|\nabla u|^2g \Big) \\
&= (1-t)\Big( \int e^{-(n+1)u} dvol_g \Big)^{\frac{2}{n+1}} +
\psi(t)f(x)e^{-2u}, \end{split}
\end{align}
where $\psi \in C^1[0,1]$ satisfies $0 \leq \psi(t) \leq 1, \psi(0)
= 0,$ and $\psi(t) \equiv 1$ for $t \leq \frac{1}{2}$; and
$\lambda_k$ is given by
\begin{align*}
\lambda_k = {n \choose k}^{-1/k}.
\end{align*}
Note that when $t=1$, equation (\ref{start}) is just equation
(\ref{hessk}).  Thus, we have constructed a deformation of
(\ref{hessk}) by connecting it to a ``less nonlinear" equation at
$t=0$:
\begin{align} \label{star0}
\sigma_k^{1/k}\big(\lambda_k g + \nabla^2 u + du \otimes du -
\frac{1}{2}|\nabla u|^2g\big) = \Big( \int e^{-(n+1)u} dvol_g
\Big)^{\frac{2}{n+1}}.
\end{align}
This 'initial' equation turns out to be much easier to analyze.
Indeed, if we assume that $g$ has been normalized to have unit
volume, then $u_0 \equiv 0$ is the unique solution of (\ref{star0}).
Since the initial equation admits a solution, one might hope to use
topological methods to establish the existence of a solution to
(\ref{start}) when $t = 1$.

In a previous paper (\cite{GVAM}) we studied (\ref{start}) with
$f(x) = constant > 0$.   Note that when $t = 0$, equation
(\ref{start}) is identical to the initial equation in \cite{GVAM};
for this reason we will only provide an outline of the degree theory
argument, and refer the reader to Section 4 of \cite{GVAM} for many
of the details.

To begin, define the operator
\begin{align} \label{psit} \begin{split}
\Psi_t[u] = \sigma_k^{1/k}& \Big( \lambda_k(1 - \psi(t))g + \psi(t)A
+
\nabla^2 u + du \otimes du - \frac{1}{2}|\nabla u|^2g \Big) \\
& - (1-t)\Big( \int e^{-(n+1)u} dvol_g \Big)^{\frac{2}{n+1}} -
\psi(t) f(x)e^{-2u}.  \end{split}
\end{align}
As we observed above, when $t = 0$
\begin{align} \label{initsoln}
\Psi_0[u_0] = 0.
\end{align}
In fact, $u_0$ is the unique solution, and the linearization of
$\Psi_0$ at $u_0$ is invertible.  It follows that the Leray-Schauder
degree $\mbox{deg}(\Psi_0,\mathcal{O}_0,0)$ defined by Li
\cite{Yanyan2} is non-zero, where $\mathcal{O}_0 \subset
C^{4,\alpha}$ is a neighborhood of the zero solution. Of course, we
would like to use the homotopy-invariance of the degree to conclude
that the degree $\mbox{deg}(\Psi_t,\mathcal{O},0)$ is non-zero for
some open set $\mathcal{O} \subset C^{4,\alpha}$; to do so we need
to establish {\it a priori} bounds for solutions of (\ref{start})
which are independent of $t$ (in order to define $\mathcal{O}$).
Using the $\epsilon$-regularity result of Guan-Wang
\cite{GuanWang1}, one can easily obtain bounds when $t < 1$:

\begin{theorem}[\cite{GVAM}, Theorem 2.1]
\label{GVAMT21}
For any fixed $0 < \delta < 1$, there is a constant $C =
C(\delta,g)$ such that any solution of (\ref{start}) with $t \in
[0,1-\delta]$ satisfies
\begin{align} \label{tless1bound}
\| u \|_{C^{4,\alpha}} \leq C.
\end{align}
\end{theorem}

The question that remains---and that we ultimately address in this
article for $k > n/2$---is the behavior of a sequence of solutions
$\{ u_{t_i} \}$ as $t_i \to 1$.  We will prove

\begin{theorem} \label{compsoln}
Let $(M^n,g)$ be a closed, compact Riemannian manifold that is not
conformally equivalent to the round $n$-dimensional sphere. Then
there is a constant $C = C(g)$ such that any solution $u$ of
(\ref{start}) satisfies
\begin{align} \label{apropt}
\| u \|_{C^{4,\alpha}} \leq C.
\end{align}
\end{theorem}

Theorem \ref{compsoln} allows us to properly define the degree of
the map $\Psi_t[\cdot]$, and by homotopy invariance we conclude the
existence of a solution of (\ref{start}) for $t = 1$.

To prove Theorem \ref{compsoln} we argue by contradiction.  Thus, we
assume $(M^n,g)$ is not conformally equivalent to the round sphere,
and let $u_i = u_{t_i}$ be a sequence of solutions to (\ref{start})
with $t_i \to 1$ and such that $\|u_i \|_{L^{\infty}} \to \infty$.
In the next section we derive various estimates used to analyze the
behavior of this sequence.

\subsection{Degree theory for other symmetric functions}
\label{degreeF}

As in the above case, we define
\begin{align} \label{psitF} \begin{split}
\tilde{\Psi}_t[u] = F& \Big( \lambda_k^{\prime}(1 - \psi(t))g +
\psi(t)A +
\nabla^2 u + du \otimes du - \frac{1}{2}|\nabla u|^2g \Big) \\
& - (1-t)\Big( \int e^{-(n+1)u} dvol_g \Big)^{\frac{2}{n+1}} -
\psi(t) f(x)e^{-2u},  \end{split}
\end{align}
where $\lambda_k^{\prime}$ is chosen so that $F \big(
\lambda_k^{\prime} \cdot g \big) = 1$.

It is very straightforward to adapt the construction above in order
to define the Leray-Schauder degree of solutions to (\ref{psitF}).
The analog of Theorem \ref{GVAMT21} is proved in a similar fashion
using the local estimates of S. Chen (Theorem \ref{SophieEstimate}).
 In the course of
proving Theorem \ref{compsoln}, in this paper we will also prove

\begin{theorem} \label{compsoln2}
Let $(M^n,g)$ be a closed, compact Riemannian manifold that is not
conformally equivalent to the round $n$-dimensional sphere. Then
there is a constant $C = C(g)$ such that any solution $u$ of
$\tilde{\Psi}_t[u] = 0$ satisfies $\| u \|_{C^{4,\alpha}} \leq C$.
\end{theorem}

\section{Local estimates}
\label{localestimates}

In this section we state some local results for solutions of
\begin{align} \label{sigmakeqn}
\sigma_k^{1/k}\big( (1-s)g + sA + \nabla^2 u + du \otimes du -
\frac{1}{2}|\nabla u|^2g \big) = \mu + f(x)e^{-2u},
\end{align}
where $u$ is assumed to be $k$-admissible, and $s \in [0,1], \mu
\geq 0$ are constants. Of course, equation (\ref{start}) is of this
form; in particular each function $\{ u_i \}$ in the sequence
defined above satisfies an equation like (\ref{sigmakeqn}).

The results of this section are of two types: the pointwise $C^1$-
and $C^2$-estimates of Guan-Wang (\cite{GuanWang1}), and various
integral estimates.  The first integral estimate (Proposition
\ref{lpmaxbound}) already appeared---albeit in a slightly different
form---in \cite{CGY1} and \cite{Gurskythesis}.

The main integral result is Theorem \ref{intgrad}.  It is a kind of
weighted $L^p$-gradient estimate that holds for $k$-admissible
metrics when $k
> n/2$, and has the advantage of assuming only minimal regularity for the metric.
This flexibility can be important when studying limits of sequences
of solutions to (\ref{sigmakeqn}), which may only be in $C^{1,1}$

\subsection{Pointwise Estimates}

Before recalling the results of Guan and Wang we should point out
that they studied equation (\ref{sigmak}) for $s = 1$ and $\mu = 0$.
However, as we explained in Section 2 of \cite{GVAM}, there is only
one line in Guan and Wang's argument that needs to be modified, and
then only slightly (see the paragraph following Lemma 2.4 in
\cite{GVAM} for details).

\begin{lemma} {\em (Theorem 1.1 of \cite{GuanWang1})}
 \label{gradrest}
Let $u \in C^4(M^n)$ be a $k-$admissible solution of
(\ref{sigmakeqn}) in $B(x_0,\rho)$, where $x_0 \in M^n$ and $\rho
> 0$.  Then there is a constant
\begin{align*}
C_1 =
C_1(k,n,\mu,\|g\|_{C^3(B(x_0,\rho))},\|f\|_{C^2(B(x_0,\rho))}),
\end{align*}
such that
\begin{align} \label{GWlocalest2}
|\nabla^2 u|(x) + |\nabla u|^2(x) \leq C_1  \big(\rho^{-2} + e^{-2
\inf_{B(x_0,\rho)}u}\big)
\end{align}
for all $x \in B(x_0,\rho/2)$.
\end{lemma}

\begin{remark}  Guan-Wang did not include the explicit dependence of
their estimates on the radius of the ball.  Since it will be
important in certain applications, we have done so here.  The
dependence is easy to establish using a typical dilation argument.
\end{remark}

An immediate corollary of this estimate is an $\epsilon-$regularity
result:

\begin{proposition} {\em (Proposition 3.6 of \cite{GuanWang1})}
\label{epsilonregr} There exist constants $\epsilon_0 > 0$ and $C =
C(g,\epsilon_0)$ such that any solution $u \in C^2(B(x_0,\rho))$ of
(\ref{sigmakeqn}) with
\begin{align}  \label{smallrenergy}
\int_{B(x_0,\rho)} e^{-nu} dvol_g \leq \epsilon_0,
\end{align}
satisfies
\begin{align} \label{lowerb}
\inf_{B(x_0,\rho/2)} u \geq -C + \log \rho.
\end{align}  Consequently, there is a constant
\begin{align*}
C_2 = C_2(k,n,\mu,\epsilon_0,\|g\|_{C^3(B(x_0,\rho))}),
\end{align*}
such that
\begin{align} \label{gradonr}
|\nabla^2 u|(x) + |\nabla u|^2(x) \leq C_2\rho^{-2}
\end{align}
for all $x \in B(x_0,\rho/4)$.
\end{proposition}

\begin{remark}  The same argument used in the proofs of Lemma \ref{gradrest} and
Proposition \ref{epsilonregr} can be used to show that any
$\Gamma$-admissible solution of (\ref{hessF}) will satisfy the
inequalities (\ref{GWlocalest2}),(\ref{lowerb}), and
(\ref{gradonr}).  Note that the homogeneity assumption on $F$ is
crucial in this respect.
\end{remark}

\subsection{Integral Estimates}

We now turn to integral estimates.  The first result shows that any
local $L^p$-bound on $e^u$ immediately gives a global $\sup$-bound.
\begin{proposition}  \label{lpmaxbound}
Let $u \in C^2(M^n)$ and assume $g_u = e^{-2u}g$ has non-negative
scalar curvature.  Suppose there is a ball $B = B(x,\rho) \subset
M^n$ and constants $\alpha_0 > 0$ and $B_0 > 0$ with
\begin{align} \label{locallp}
\int_{B(x,\rho)} e^{\alpha_0 u} dvol_g \leq B_0.
\end{align}
Then there is a constant
\begin{align*}
C = C(g,vol(B(x,\rho)),\alpha_0,B_0),
\end{align*}
such that
\begin{align} \label{maxfromlp}
\max_{M^n} u \leq C.
\end{align}
\end{proposition}
\begin{proof}
If $R_u$ denotes the scalar curvature of $g_u$, then (\ref{Achange})
implies
\begin{align*}
\frac{1}{2(n-1)}R_ue^{-2u} = \frac{1}{2(n-1)}R + \Delta u -
\frac{(n-2)}{2}|\nabla u|^2.
\end{align*}
Therefore, if $R_u \geq 0$ we conclude
\begin{align} \label{Scaleqn}
-\Delta u + \frac{(n-2)}{2}|\nabla u|^2 + \frac{1}{2(n-1)}R \geq 0,
\end{align}
hence
\begin{align} \label{Rnonneg}
-\Delta u + \frac{(n-2)}{2}|\nabla u|^2 \geq -C.
\end{align}
In what follows, it will simplify our calculations if we let $v =
e^{-\frac{(n-2)}{2}u}$.  In terms of $v$, the bound (\ref{locallp})
becomes
\begin{align} \label{locallpv}
\int_{B(x,\rho)} v^{-p_0} dvol_g \leq B_0
\end{align}
where $p_0 = \frac{(n-2)}{2}\alpha_0$.  Also, inequality
(\ref{Scaleqn}) becomes
\begin{align} \label{Scaleqnv}
\Delta v - \frac{(n-2)}{4(n-1)}Rv \leq 0,
\end{align}
and (\ref{Rnonneg}) becomes
\begin{align} \label{Rlin}
\Delta v \leq Cv.
\end{align}
It follows that any {\it global} $L^p$-bound of the form
(\ref{locallpv}) implies a lower bound on $v$ (and therefore an
upper bound on $u$). That is, if $p > 0$, then
\begin{align} \label{lptoinf}
\int_{M^n} v^{-p} dvol_g \leq C_0 \Rightarrow \inf_{M^n} v \geq
C(C_0) > 0.
\end{align}
There are various ways to see this; for example, by using the
Green's representation. It therefore remains to prove that one can
pass from the local $L^p$-bound of (\ref{locallpv}) to a global one:
\begin{lemma} \label{lploctoglobe}
For $p \in (0,p_0)$ sufficiently small, there is a constant
\begin{align*}
C = C(g, vol(B(x,\rho)), p, B_0),
\end{align*}
such that
\begin{align} \label{globlp}
\int_{M^n} v^{-2p} dvol_g \leq C.
\end{align}
\end{lemma}
\begin{proof}  This is Lemma 4.3 of \cite{CGY1} (see also
Lemma 4.1 of \cite{Gurskythesis}).
\end{proof}
This completes the proof of Proposition \ref{lpmaxbound}.
\end{proof}
\begin{remark} If $u \in C^2$ is a $\Gamma$-admissible solution of
(\ref{hessF}), then by definition the scalar curvature of $g_u =
e^{-2u}g$ is positive.  Therefore, Proposition \ref{lpmaxbound} is
applicable.
\end{remark}

The next result is an integral gradient estimate for admissible
metrics.  Before we give the precise statement, a brief remark is
needed about the regularity assumptions of the result and their
relationship to curvature.

If $u \in C^{1,1}$, then Rademacher's Theorem says that the Hessian
of $u$ is defined almost everywhere, and therefore by
(\ref{Achange}) the Schouten tensor $A_u$ of $g_u$ is defined almost
everywhere.  In particular, the notion of $k$-admissibility
(respectively, $\Gamma$-admissibility) can still be defined: it
requires that the eigenvalues of $A_u$ are in
$\Gamma_k^{+}(\mathbf{R}^n)$ ({\it resp.}, $\Gamma$) at almost every
$x \in M^n$. Likewise, the condition of non-negative Ricci curvature
(a.e.) is well defined.
\begin{theorem} \label{intgrad}  Let $u \in
C_{loc}^{1,1}\big(A(\frac{1}{2}r_1,2r_2)\big)$, where $x_0 \in M^n$
and $A(\frac{1}{2}r_1,2r_2)$ denotes the annulus
$A(\frac{1}{2}r_1,2r_2) \equiv B(x_0,2r_2) \setminus
\overline{B(x_0, \frac{1}{2}r_1)}$, with $0 < r_1 < r_2$. Assume
$g_u = e^{-2u}g$ satisfies
\begin{align} \label{posRicci}
Ric(g_u) - 2\delta \sigma_1(A_u)g \geq 0
\end{align}
almost everywhere in $A(\frac{1}{2}r_1,2r_2)$ for some $0 \leq
\delta < \frac{1}{2}$. Define
\begin{align} \label{alphadeltadef}
\alpha_{\delta} = \frac{(n-2)}{(1-2\delta)}\delta \geq 0.
\end{align}
Then given any $\alpha > \alpha_{\delta}$, there are constants $p
\geq n$ and $C = C((\alpha - \alpha_{\delta})^{-1},n) > 0$ such that
\begin{align} \label{localintgrad0} \begin{split}
\int_{A(r_1,r_2)} |\nabla u|^p e^{\alpha u} & dvol_g \leq C\Bigg(
\int_{A(\frac{1}{2}r_1,2r_2)} |Ric_g|^{p/2} e^{\alpha u} dvol_g \\
&+ r_1^{-p} \int_{A(\frac{1}{2}r_1,r_1)}e^{\alpha u}dvol_g +
r_2^{-p} \int_{A(r_2, 2r_2)} e^{\alpha u} dvol_g \Bigg).
\end{split}
\end{align}
In fact, we can take
\begin{align} \label{pdef}
p = n + 2\alpha_{\delta} \geq n.
\end{align}
\end{theorem}
Now suppose $g_u = e^{-2u}g$ is $k$-admissible, with $k > n/2$. By a
result of Guan-Viaclovsky-Wang (\cite{GVW}, Theorem 1), inequality
(\ref{posRicci}) holds for any $\delta$ satisfying $\delta \leq
\frac{(2k -n)(n-1)}{2n(k-1)}$.  Therefore,
\begin{corollary}
\label{intgradlocalcor} Let $u \in
C_{loc}^{1,1}\big(A(\frac{1}{2}r_1,2r_2)\big)$, where $x_0 \in M^n$
and $A(\frac{1}{2}r_1,2r_2)$ denotes the annulus
$A(\frac{1}{2}r_1,2r_2) \equiv B(x_0,2r_2) \setminus \overline{
B(x_0, \frac{1}{2}r_1)}$, with $0 < r_1 < r_2$. Assume $g_u =
e^{-2u}g$ is $k$-admissible with $k > n/2$.  Suppose $\delta \geq 0$
satisfies
\begin{align*}
0 \leq & \delta \leq \min \Big\{ \frac{1}{2}, \frac{(2k -
n)(n-1)}{2n(k-1)}\Big\},
\end{align*}
and define
\begin{align*}
\alpha_{\delta} = \frac{(n-2)}{(1-2\delta)}\delta.
\end{align*}
Then given any $\alpha > \alpha_{\delta}$, there are constants $p
> n$ and $C = C((\alpha - \alpha_{\delta})^{-1},n) > 0$ such that
\begin{align} \label{localintgrad}  \begin{split}
\int_{A(r_1,r_2)} |\nabla u|^p e^{\alpha u} & dvol_g \leq C\Bigg(
\int_{A(\frac{1}{2}r_1,2r_2)} |Ric_g|^{p/2} e^{\alpha u} dvol_g \\
&+ r_1^{-p} \int_{A(\frac{1}{2}r_1,r_1)}e^{\alpha u}dvol_g +
r_2^{-p} \int_{A(r_2, 2r_2)} e^{\alpha u} dvol_g \Bigg).
\end{split}
\end{align}
\end{corollary}
We now give the proof of Theorem \ref{intgrad}.
\begin{proof}
Let
\begin{align} \label{Sdef}
S = S_g =  Ric_g - 2\delta\sigma_1(A_g)g.
\end{align}
By the curvature transformation formula (\ref{Achange}), for any
conformal metric $g_u = e^{-2u}g$ the relationship between $S_u =
S_{g_u}$ and $S_g$ is given by
\begin{align} \label{Schange}
S_u = (n-2)\nabla^2 u + (1-2\delta)\Delta u g + (n-2)du \otimes du -
(n-2)(1 - \delta)|\nabla u|^2 g + S_g.
\end{align}
Now assume $g_u$ is a metric for which $S_u \geq 0$ a.e. in $A =
A(\frac{1}{2}r_1,2r_2)$.  From (\ref{Schange}), it follows that
\begin{align} \label{gradfact} \begin{split}
0 &\leq \frac{1}{(n-2)}S_u(\nabla u,\nabla u) \\
&= \nabla^2 u (\nabla u, \nabla u) + \frac{(1-2\delta)}{(n-2)}\Delta
u |\nabla u|^2 + \delta |\nabla u|^4 + \frac{1}{(n-2)}S_g(\nabla
u,\nabla u)
\end{split}
\end{align}
a.e. in $A$.  Using this identity we have
\begin{lemma}
\label{DivForm} For any $\alpha \in \mathbf{R}, \delta \geq 0$, $u$
satisfies
\begin{align} \label{Divu}
\nabla_i (|\nabla u|^{p-2}e^{\alpha u} \nabla_i u) \geq (\alpha -
\alpha_{\delta})|\nabla u|^p e^{\alpha u} -
\frac{1}{(1-2\delta)}|S_g||\nabla u|^{p-2}e^{\alpha u}
\end{align}
almost everywhere in $A$, where $\alpha_{\delta}$ and $p$ are
defined in (\ref{alphadeltadef}) and (\ref{pdef}), respectively.
\end{lemma}
\begin{proof}  Since $u \in C_{loc}^{1,1}$ and $p > 2$, the vector field
\begin{align} \label{Xdef}
X = |\nabla u|^{p-2}e^{\alpha u} \nabla u
\end{align}
is locally Lipschitz. Therefore, its divergence is defined at any
point where the Hessian of $u$ is defined--in particular, almost
everywhere in $A$. Moreover, at any point where (\ref{gradfact}) is
valid we have
\begin{align} \label{divu1} \begin{split}
\nabla_i (|\nabla u|^{p-2}e^{\alpha u} \nabla_i u) &= \nabla_i
(|\nabla u|^{p-2}) e^{\alpha u} \nabla_i u + |\nabla
u|^{p-2}\nabla_i(e^{\alpha u})\nabla_i u + |\nabla
u|^{p-2}e^{\alpha u} \Delta u \\
&= (p-2)|\nabla u|^{p-4}e^{\alpha u} \nabla^2 u (\nabla u,\nabla u)
+ |\nabla u|^{p-2}e^{\alpha u} \Delta u + \alpha |\nabla u|^p
e^{\alpha u}.
\end{split}
\end{align}
Since $p > 2$ we can substitute inequality (\ref{gradfact}) into
(\ref{divu1}) to get
\begin{align} \label{divu2} \begin{split}
\nabla_i (|\nabla u|^{p-2}& e^{\alpha u} \nabla_i u) \geq \big[1
-(p-2)\frac{(1-2\delta)}{(n-2)} \big]|\nabla u|^{p-2}e^{\alpha u}
\Delta u \\
& + \big[ \alpha - \delta(p-2)\big] |\nabla u|^p e^{\alpha u} -
\frac{(p-2)}{(n-2)}|S_g||\nabla u|^{p-2} e^{\alpha u}.
\end{split}
\end{align}
Then (\ref{Divu}) follows from (\ref{divu2}) by taking
\begin{align*}
p = \frac{n - 4\delta}{1-2\delta} = n + 2\alpha_{\delta}.
\end{align*}
\end{proof}
Let $\eta$ be a cut-off function satisfying
\begin{align*}
\eta(x) = \left\{ \begin{array}{lll}
0 & x \in B(x_0,\frac{3}{4}r_1), \\
1 & x \in A(r_1,r_2), \\
0 & x \in B(x_0,2r_2) \setminus B(x_0,\frac{3}{2}r_2),
\end{array}
\right.
\end{align*}
and
\begin{align*}
|\nabla \eta(x)| \leq \left\{ \begin{array}{ll}
Cr_1^{-1} & x \in A(\frac{1}{2}r_1,r_1), \\
Cr_2^{-1} & x \in A(r_2,2r_2),
\end{array}
\right.
\end{align*}
and $|\nabla \eta| = 0$ otherwise.  Multiplying both sides of
(\ref{Divu}) by $\eta^p$ and applying the divergence theorem (which
is valid since the vector field $X$ in (\ref{Xdef}) is Lipschitz),
we get
\begin{align} \label{inteq1} \begin{split}
\int -p \langle \nabla \eta, \nabla u \rangle & \eta^{p-1} |\nabla
u|^{p-2}e^{\alpha u} dvol_g \geq \\ &(\alpha - \alpha_{\delta}) \int
|\nabla u|^p e^{\alpha u} \eta^p dvol_g - \frac{1}{(1-2\delta)} \int
|S_g||\nabla u|^{p-2}e^{\alpha u} \eta^p dvol_g. \end{split}
\end{align}
Using the obvious bound $|S_g| \leq C|Ric_g|$ and rearranging terms
in (\ref{inteq1}) gives
\begin{align} \label{inteq2} \begin{split}
 \int & |\nabla u|^p  e^{\alpha u} \eta^p
dvol_g \leq \\ & C\big((\alpha -
\alpha_{\delta})^{-1},\delta,n\big)\Big[ \int |Ric_g| |\nabla
u|^{p-2}e^{\alpha u}\eta^p dvol_g + \int |\nabla \eta| |\nabla
u|^{p-1} e^{\alpha u} \eta^{p-1} dvol_g \Big].  \end{split}
\end{align}
Applying H\"older's inequality to the two integrals on the $RHS$ we
have
\begin{align*}
\int |\nabla u|^p & e^{\alpha u} \eta^p dvol_g \\ &\leq
C(\theta,\delta,n) \Bigg[  \Bigg( \int |\nabla u|^p e^{\alpha u}
\eta^p dvol_g \Bigg)^{(p-2)/p} \Bigg( \int |Ric_g|^{p/2} e^{\alpha
u} \eta^p dvol_g \Bigg)^{2/p} \\ &\quad + \Bigg( \int |\nabla u|^p
e^{\alpha u} \eta^p dvol_g \Bigg)^{(p-1)/p}\Bigg( \int |\nabla
\eta|^p e^{\alpha u} dvol_g \Bigg)^{1/p} \Bigg],
\end{align*}
which implies
\begin{align*}
\int |\nabla u|^p e^{\alpha u} \eta^p dvol_g  &\leq  C\Big[
  \int |Ric_g|^{p/2} e^{\alpha u} \eta^p dvol_g \\
&\quad +  \int |\nabla \eta|^p e^{\alpha u} dvol_g \Big].
\end{align*}
Therefore, using the properties of the cut-off function $\eta$ we
conclude
\begin{align*}
\int_{A(r_1,r_2)} |\nabla u|^p e^{\alpha u} dvol_g &\leq
\int |\nabla u|^p e^{\alpha u} \eta^p dvol_g \\
&\leq C\Big[
  \int |Ric_g|^{p/2} e^{\alpha u} \eta^p dvol_g +  \int |\nabla \eta|^p e^{\alpha u} dvol_g \Big] \\
&\leq C \int_{A(\frac{1}{2}r_1,2r_2)} |Ric_g|^{p/2} e^{\alpha u}
dvol_g \\
& \quad + C r_1^{-p} \int_{A(\frac{1}{2}r_1,r_1)} e^{\alpha u}
dvol_g + C r_2^{-p} \int_{A(r_2,2r_2)} e^{\alpha u} dvol_g.
\end{align*}
This completes the proof of Theorem \ref{intgrad}.
\end{proof}
Theorem \ref{intgrad} has a global version:
\begin{corollary} \label{intgradglobalcor} Let $u \in
C_{loc}^{1,1}(M^n)$, and assume $g_u = e^{-2u}g$ is $k$-admissible
with $k > n/2$.  Suppose $\delta \geq 0$ satisfies
\begin{align*}
0 \leq & \delta \leq \min \Big\{ \frac{1}{2}, \frac{(2k -
n)(n-1)}{2n(k-1)}\Big\},
\end{align*}
and define
\begin{align*}
\alpha_{\delta} = \frac{(n-2)}{(1-2\delta)}\delta.
\end{align*}
Then given any $\alpha > \alpha_{\delta}$, there are constants $p
\geq n$ and $C = C((\alpha - \alpha_{\delta})^{-1},n,g) > 0$ such
that
\begin{align} \label{globalintgrad}
\int_{M^n} |\nabla u|^p e^{\alpha u}dvol_g \leq C
\int_{M^n}e^{\alpha u}dvol_g.
\end{align}
\end{corollary}
Using the Sobolev Imbedding Theorem, Corollary
\ref{intgradglobalcor} implies a pointwise estimate:
\begin{corollary}
\label{holdercor} Let $u \in C_{loc}^{1,1}(M^n)$, and assume $g_u =
e^{-2u}g$ is $k$-admissible with $k > n/2$.  Suppose $\delta > 0$
satisfies
\begin{align} \label{deltaineq}
0 <  \delta \leq \min \Big\{ \frac{1}{2}, \frac{(2k -
n)(n-1)}{2n(k-1)}\Big\},
\end{align}
and define
\begin{align}
\alpha_{\delta} = \frac{(n-2)}{(1-2\delta)}\delta.
\end{align}
Then given any $\alpha > \alpha_{\delta}$, there is a constant $C =
C((\alpha - \alpha_{\delta})^{-1},n,g) > 0$ such that
\begin{align} \label{lptoholder}
\| e^{(\alpha/p)u} \|_{C^{\gamma_0}} \leq C \| e^{(\alpha/p)u}
\|_{L^p},
\end{align}
where
\begin{align} \label{gammadef}
\gamma_0 =  \ds \frac{2\alpha_{\delta}}{n + 2 \alpha_{\delta} } > 0.
\end{align}
\end{corollary}
\begin{proof}
Inequality (\ref{globalintgrad}) implies
\begin{align*}
\| e^{(\alpha/p)u} \|_{W^{1,p}} \leq C \| e^{(\alpha/p)u} \|_{L^p},
\end{align*}
where $W^{1,p}$ denotes the Sobolev space of functions $\varphi \in
L^p$ with $|\nabla \varphi| \in L^p$.  The Sobolev Imbedding Theorem
implies the bound (\ref{lptoholder}) for
\begin{align*}
\gamma_0 = 1 - \frac{n}{p}.
\end{align*}
If we take $p = n + 2\alpha_{\delta}$ as in (\ref{pdef}), then
(\ref{gammadef}) follows.
\end{proof}

\subsection{Local estimates for other symmetric functions}

As we observed in the Remarks following the proofs of Propositions
\ref{epsilonregr} and \ref{lpmaxbound}, any $\Gamma$-admissible
solution of (\ref{hessF}) automatically satisfies the conclusions of
Lemma \ref{gradrest} and Propositions \ref{epsilonregr} and
\ref{lpmaxbound}. Furthermore, the condition (\ref{Ricposcond})
implies that any $\Gamma$-admissible solution satisfies inequality
(\ref{posRicci}) of Theorem \ref{intgrad}.  Therefore, this result
and its corollaries remain valid for $\Gamma$-admissible solutions.

Furthermore, suppose $\{ u_i \}$ is a sequence of solutions to
$\tilde{\Psi}_{t_i}[u_i] = 0$, as described in Section
\ref{degreeF}.  Then the conclusion of Proposition \ref{floordrops}
also holds for this sequence, since the proof just relies on the
local estimates of S. Chen.

\section{The blow-up}
\label{blow-up}

In this section we begin a careful analysis of a sequence $\{ u_i
\}$ of solutions to (\ref{start}). We may assume that $u_{t_i} =
u_i$ with $t_i
> 1/2$; this implies $\psi(t_i) = 1$, so (\ref{start}) becomes
\begin{align} \label{start1} \begin{split}
\sigma_k^{1/k}& \Big( A +
\nabla^2 u_i + du_i \otimes du_i - \frac{1}{2}|\nabla u_i|^2g \Big) \\
&= (1-t_i)\Big( \int e^{-(n+1)u_i} dvol_g \Big)^{\frac{2}{n+1}} +
f(x) e^{-2u_i}. \end{split}
\end{align}
In particular, the metrics $g_i = e^{-2u_i}g$ are $k$-admissible.

Our first observation is that the sequence $\{ u_i \}$ must have
$\min u_i \to -\infty$.

\begin{proposition} \label{floordrops}
If there is a lower bound
\begin{align}  \label{floor}
\min u_i \geq -C
\end{align}
then there is an upper bound as well; i.e.,
\begin{align*}
\| u_i \|_{L^{\infty}} \leq C.
\end{align*}
Therefore, we must have
\begin{align} \label{infdown}
\min u_i \to -\infty.
\end{align}
\end{proposition}
\begin{proof}
From \cite[Lemma 2.5]{GVAM} we know that $\min u_i$ is bounded
above. Therefore, assuming (\ref{floor}) holds, we have
\begin{align} \label{bothsides}
-C \leq \min u_i \leq C.
\end{align}
By Lemma \ref{gradrest}, the lower bound (\ref{floor}) implies a
gradient bound, which gives a Harnack inequality of the form
\begin{align} \label{littleH}
\max u_i \leq \min u_i + C'.
\end{align}
Combining (\ref{bothsides}) and (\ref{littleH}) we conclude
\begin{align*}
\max u_i \leq C + C'.
\end{align*}
Consequently, (\ref{infdown}) must hold.
\end{proof}
The next lemma will be used to show that the sequence $\{ u_i \}$
can only concentrate at finitely many points:

\begin{lemma} \label{ricvolume}
The volume and Ricci curvature of the metrics $\{ g_i \}$ satisfy
\begin{align} \label{volumeabove}
vol(g_i) \leq v_0,
\end{align}
\begin{align} \label{Ricnonneg}
Ric_{g_i} \geq c_0 g_i,
\end{align}
where $c_0 > 0$.
\end{lemma}

\begin{proof} Since each $g_i$ is $k$-admissible with $k > n/2$,
by Theorem 1 of \cite{GVW} the Ricci curvature of $g_i$ satisfies
\begin{align} \label{GVWRiclower}
Ric_{g_i} \geq \frac{(2k-n)(n-1)}{(k-1)}{n \choose
k}^{-1/k}\sigma_k^{1/k}(A_{g_i})g_i.
\end{align}
Using equation (\ref{start1}) this implies
\begin{align} \label{Ricilower} \begin{split}
Ric_{g_i} &\geq \frac{(2k-n)(n-1)}{(k-1)}{n \choose
k}^{-1/k}f(x)g_i \\
&\geq c(k,n,\min f)g_i  \end{split}
\end{align}
for some $c_0(k,n,\min f) > 0$. This proves inequality
(\ref{Ricnonneg}). Also, by Bishop's Theorem (\cite{Bishop}) it
gives an upper bound on the volume of $g_i$, proving
(\ref{volumeabove}).

\end{proof}

\begin{remark} It is easy to see that inequalities (\ref{volumeabove}) and
(\ref{Ricnonneg}) are valid for any sequence $\{ u_i \}$ of
$\Gamma$-admissible solutions to (\ref{hessF}). This follows from
(\ref{Ricposcond}), (\ref{Ricposintro}), and a lower bound for the
scalar curvature.  More precisely, we claim that
\begin{align}  \label{lmi}
F(\lambda) \leq C_0 \sigma_1(\lambda)
\end{align}
for some constant $C_0 > 0$ and any $\lambda \in \Gamma$.  To see
this, by the homogeneity of $F$ it suffices to prove that it holds
for $\lambda \in \hat{\Gamma} = \{ \lambda \in \overline{\Gamma} :
|\lambda| = 1 \}$.  Now, (\ref{Ricposcond}) implies that
$\sigma_1(\lambda)
> 0$ for $\lambda \in \hat{\Gamma}$, because if $\sigma_1(\hat{\lambda}) = 0$ for some
$\hat{\lambda} \in \hat{\Gamma}$ then by (\ref{Ricposcond}) we would
have $\hat{\lambda} = 0$, a contradiction. Therefore,
$\sigma_1(\lambda) \geq c_0 > 0$ on $\hat{\Gamma}$, so
$F(\lambda)/\sigma_1(A)$ is bounded above on $\hat{\Gamma}$, which
proves (\ref{lmi}). Appealing to (\ref{Ricposintro}) and
(\ref{hessF}), we see that the Ricci curvature of $g_i = e^{-2u_i}g$
satisfies
\begin{align*}
Ric_{g_i} \geq 2 C_0^{-1}\delta f(x) g_i \geq c g_i,
\end{align*}
for some $c > 0$. Arguing as we did above, we obtain
(\ref{volumeabove}) and (\ref{Ricnonneg}).
\end{remark}
\vskip.25in

Given $x \in M^n$, define the {\it mass} of $x$ by
\begin{align} \label{massdef}
m (x) = m(\{u_i\};x) = \lim_{r \to 0} \limsup_{i\to\infty}
\int_{B(x,r)} e^{-nu_i} dvol_g.
\end{align}
This notation is intended to emphasize the dependence of the mass on
the sequence $\{ u_i \}$.  In particular, if we restrict to a
subsequence (as we will soon do), the mass of a given point may
decrease.

The $\epsilon$-regularity result of Proposition \ref{epsilonregr}
implies that, on a subsequence, only finitely many points may have
non-zero mass:

\begin{proposition}[\cite{Gurskythesis}, Section 2]
\label{singsetsize} The set $\Sigma[\{ u_i \}] = \{ x \in M^n \vert
m(x) \neq 0 \}$ is non-empty. In addition, there is a subsequence
(still denoted $\{ u_i \}$) such that with respect to this
subsequence $\Sigma$ is non-empty and consists of finitely many
points: $\Sigma = \Sigma[\{ u_i \}] = \{ x_1, x_2, \dots , x_{\ell}
\}$.
\end{proposition}

\subsection{Behavior away from the singular set}

While the sequence $\{ u_i \}$ is concentrating at the points $\{
x_1, x_2, \dots ,x_{\ell} \}$, away from these points the $u_i$'s
remain bounded from below, and the derivatives up to order two are
uniformly bounded:

\begin{proposition} \label{nearsigv1}
Given compact $K \subset M^n \setminus \Sigma$, there is a constant
$C = C(K)>0$ that is independent of $i$ such that
\begin{align} \label{infoffsigv1}
\min_{K} u_i \geq -C(K),
\end{align}
\begin{align} \label{gradoffsigv1}
\max_{K} \big[ |\nabla^2 u_i| + |\nabla u_i|^2 \big] \leq C(K),
\end{align}
for all $i \geq J = J(K)$.
\end{proposition}

\begin{proof}  Given $x \in K$, since $m(x) = 0$ there is radius
$r = r_x$ such that
\begin{align*}
\int_{B(x,r_x)} e^{-nu_i} dvol_g \leq \frac{1}{2}\epsilon_0
\end{align*}
for all $i \geq J = J_x$.  By Proposition \ref{epsilonregr} there is
a constant $C = C(r_x) > 0$ such that
\begin{align} \label{infonballv1}
\inf_{B(x,r_x/2)} u_i \geq -C,
\end{align}
\begin{align} \label{gradonballv1}
\sup_{B(x,r_x/2)} \big[ |\nabla^2 u_i| + |\nabla u_i|^2 \big] \leq
C,
\end{align}
for all $i \geq J_x$.  The balls $\{ B(x,r_x/2) \}_{x \in K}$ define
an open cover of $K$, and since $K$ is compact we can extract a
finite subcover $K \subset \cup_{\nu = 1}^{N} B(x_{\nu},r_{\nu}/2)$.
Let $J = \max_{1 \leq \nu \leq N} J_{\nu}$; then (\ref{infonballv1})
and (\ref{gradonballv1}) imply that (\ref{infoffsigv1}) and
(\ref{gradoffsigv1}) hold for all $i \geq J$.
\end{proof}

Now, fix a ``regular point" $x_0 \notin \Sigma$.  There are two
possibilities to consider, depending on whether
\begin{align} \label{I}
\limsup_{i} u_i(x_0) < +\infty
\end{align}
or
\begin{align} \label{II}
\limsup_{i} u_i(x_0) = +\infty
\end{align}
(recall that (\ref{infoffsigv1}) only provides a {\it lower} bound
for the sequence off the singular set).   These possibilities
reflect different scenarios for the convergence of (a subsequence
of) $\{ u_i \}$ on $M^n \setminus \Sigma$.   If (\ref{I}) holds, it
will be possible to extract a subsequence that converges on compact
subsets of $M^n \setminus \Sigma$ to a smooth limit $u \in
C^{\infty}(M^n \setminus \Sigma)$.  But if (\ref{II}) holds, a
subsequence diverges to $+\infty$ uniformly on compact subsets of
$M^n \setminus \Sigma$.  As we shall see, the integral gradient
estimate (Corollary \ref{holdercor}) can be used to preclude
(\ref{I}).

To this end, assume
\begin{align} \label{badcase}
\limsup_{i} u_i(x_0) < +\infty.
\end{align}
Then if $K \subset M^n \setminus \Sigma$ is a compact set containing
$x_0$, the bounds (\ref{infoffsigv1}), (\ref{gradoffsigv1}), and
(\ref{badcase}) imply there is a constant $C = C(K)
> 0$ such that
\begin{align} \label{supisbounded}
\max_{K} \Big[ |\nabla^2 u_i| + |\nabla u_i|^2 + |u_i| \Big] \leq
C(K)
\end{align}
for all $i \geq J = J(K)$.   This estimate implies that equation
(\ref{start1}) is {\it uniformly elliptic on} $K$.  Since it is
concave, by the results of Evans (\cite{Evans}) and Krylov
(\cite{Krylov}) one obtains interior $C^{2,\gamma}$-bounds for
solutions.  The Schauder interior estimates then give estimates on
derivatives of all orders:  More precisely, given $K^{\prime}
\subset K$, $m \geq 1$, and $\gamma \in (0,1)$, there is a constant
$C = C(K^{\prime}, m, \gamma)$ such that
\begin{align} \label{localholder}
\| u_i \|_{C^{m,\gamma}(K^{\prime})} \leq C.
\end{align}
After applying a standard diagonal argument, we may extract a
subsequence $u_i \to u \in C^{\infty}(M^n \setminus \Sigma)$, where
the convergence is in $C^m$ on compact sets away from $\Sigma$.

As we observed above, when restricting to subsequences it is
possible that one reduces the singular set.  However, it is always
possible to choose a subsequence of $\{ u_i \}$ and a sequence of
points $\{ P_i \}$ with
\begin{align} \label{badpoint}
\lim_{i} u_i(P_i) = -\infty, \quad P_i \to P \in \Sigma,
\end{align}
say $P = x_1$.  That is, we can always choose a subsequence for
which there is at least one singular point.  For, if such a choice
were impossible, then the original sequence would have a uniform
lower bound, and this would violate the conclusions of Proposition
\ref{floordrops}.

Using Proposition \ref{lpmaxbound} and Corollary \ref{holdercor}, we
can obtain more precise information on the behavior of the limit $u$
near the singular point $x_1$.

\begin{proposition} \label{lifenearx1}  Under the assumption
(\ref{badcase}), the function $u = \lim_{i} u_i$ has the following
properties:

\vskip.1in \noindent$(i)$  There is a constant $C_1 > 0$ such that
\begin{align} \label{uaboveisok}
\sup_{M^ \setminus \Sigma} u \leq C_1.
\end{align}

\vskip.1in \noindent $(ii)$  There is a neighborhood $U$ containing
$x_1$ with the following property:  Given $\theta
> 0$, there is a constant $C = C(\theta)$ such that
\begin{align} \label{lograte}
u(x) \leq (2 - \theta)\log d_g(x,x_1) + C(\theta)
\end{align}
for all $x \neq x_1$ in $U$.
\end{proposition}

\begin{proof}
To prove (\ref{uaboveisok}), let $K \subset M^n \setminus \Sigma$ be
a compact set containing $x_0$.  By (\ref{supisbounded}), we have a
bound
\begin{align*}
\int_{K} e^{\alpha u_i} dvol_g \leq C_{\alpha}
\end{align*}
for any $\alpha > 0$.  Therefore, by Proposition \ref{lpmaxbound} we
have a global bound
\begin{align} \label{uibound}
\max_{M^n} u_i \leq C_1.
\end{align}
Thus the limit must satisfy
\begin{align*}
\sup_{M^n \setminus \Sigma} u \leq C_1.
\end{align*}

Turning to the proof of (\ref{lograte}), note the bound
(\ref{uibound}) allows us to apply Corollary \ref{holdercor}.
Therefore, for fixed $\delta > 0$ satisfying (\ref{deltaineq}) and
any $\alpha > \alpha_{\delta} \equiv
\frac{(n-2)}{(1-2\delta)}\delta$, we have
\begin{align} \label{holderforu}
\| e^{(\alpha/p)u_i}\|_{C^{\gamma_0}} \leq C \| e^{(\alpha/p)u_i}
\|_{L^p} \leq C,
\end{align}
where
\begin{align*}
p = n + 2\alpha_{\delta},
\end{align*}
\begin{align*}
\gamma_0 =  \ds \frac{2\alpha_{\delta}}{n + 2 \alpha_{\delta} } > 0.
\end{align*}
Choose a small neighborhood $U$ of $x_1$ that is disjoint from the
other singular points.   For $i > J$ sufficiently large we may
assume that $P_i \in U$, where $P_i \to x_1$ is the sequence in
(\ref{badpoint}).  If $x \neq x_1$ is a point in $U$, then by
(\ref{holderforu})
\begin{align} \label{holdvsdist}
| e^{(\alpha/p) u_i(x)} - e^{(\alpha/p) u_i(P_i)} | \leq C
d_g(x,P_i)^{\gamma_0}.
\end{align}
Letting $i \to \infty$ in (\ref{holdvsdist}), by (\ref{badpoint}) we
conclude
\begin{align*}
e^{ (\alpha/p)u(x) } \leq C d_g(x,x_1)^{\gamma_0}.
\end{align*}
This implies
\begin{align} \label{rate1}
e^{u(x)} \leq C d_g(x,x_1)^{p\gamma_0/\alpha}.
\end{align}
Using the definition of $p$ in (\ref{pdef}), the exponent in
(\ref{rate1}) satisfies
\begin{align*}
\ds \frac{p\gamma_0}{\alpha} = 2\frac{\alpha_{\delta}}{\alpha} < 2,
\end{align*}
so taking logarithms in (\ref{rate1}) we get
\begin{align} \label{lograte2}
u(x) \leq \ds 2\big(\frac{\alpha_{\delta}}{\alpha}\big) \log
d_g(x,x_1) + C.
\end{align}
Therefore, given $\theta > 0$, we can choose $\alpha >
\alpha_{\delta}$ close enough to $\alpha_{\delta}$ so that
\begin{align*}
2\frac{\alpha_{\delta}}{\alpha} \geq 2 - \theta,
\end{align*}
and (\ref{lograte}) follows from (\ref{lograte2}).

\end{proof}

While Proposition \ref{lifenearx1} gives fairly precise {\it upper}
bounds on $u = \lim_{i} u_i$ near $x_1$, the epsilon-regularity
result Proposition \ref{epsilonregr} can be used to give {\it lower}
bounds:

\begin{proposition} \label{lifenearx1frombelow}
There is a neighborhood $U^{\prime}$ of $x_1$ and a constant $C > 0$
such that
\begin{align} \label{eregforu}
u(x) \geq \log d_g(x,x_1) - C
\end{align}
for all $x \neq x_1$ in $U^{\prime}$.
\end{proposition}

\begin{proof}  This result follows from inequality (\ref{lowerb}).
More precisely, by the volume bound (\ref{volumeabove}),
\begin{align*}
vol(g_i) = \int_{M^n} e^{-nu_i} dvol_g \leq v_0.
\end{align*}
It follows that $u = \lim_{i} u_i$ satisfies
\begin{align*}
vol(g_u) = \int_{M^n} e^{-nu} dvol_g \leq v_0.
\end{align*}
Therefore, for $\rho_0 > 0$ small enough,
\begin{align*}
\int_{B(x_1,\rho_0)} e^{-nu} dvol_g \leq \frac{1}{2}\epsilon_0,
\end{align*}
where $\epsilon_0$ is the constant in the statement of Proposition
\ref{epsilonregr}.

Given $x \neq x_1$ in $U^{\prime} = B(x_1,\rho_0/2)$, let $\rho =
\frac{1}{2}d_g(x,x_1)$.  Then $B(x,\rho) \subset B(x_1,\rho_0)$, so
\begin{align*}
\int_{B(x,\rho)} e^{-nu} dvol_g < \int_{B(x_1,\rho_0)} e^{-nu}
dvol_g \leq \frac{1}{2}\epsilon_0.
\end{align*}
Therefore, by inequality (\ref{lowerb}),
\begin{align*}
\inf_{B(x,\rho/2)} u \geq \log \rho - C,
\end{align*}
which implies
\begin{align*}
u(x) \geq \log d_g(x,x_1) - C.
\end{align*}
\end{proof}

Combining Propositions \ref{lifenearx1} and
\ref{lifenearx1frombelow}, we conclude that the assumption
(\ref{badcase}) on which they are based can not be true:

\begin{corollary} \label{uitoplusinfty}
The sequence $\{ u_i \}$ must satisfy
\begin{align} \label{supatpointbig}
\limsup_{i} u_i(x_0) = +\infty.
\end{align}
Consequently, there is a subsequence (again denoted $\{ u_i \}$) and
a non-empty set of points $\Sigma_0 = \{ x_1, x_2, \dots, x_{\nu} \}
\subseteq \Sigma$ with the following properties:

\vskip.1in \noindent $(i)$  Given any compact $K \subset M^n
\setminus \Sigma_0$ and number $N
>0$, there is a $J = J(K,N)$ such that
\begin{align} \label{infisbig}
\min_{K} u_i \geq N
\end{align}
for all $i \geq J = J(K,N)$.

\vskip.1in \noindent $(ii)$  For each $x_k \in \Sigma_0$, there is a
sequence of point $\{ P_{k,i} \}$ such that
\begin{align} \label{badpoints}
\lim_i P_{k,i} = x_k,
\end{align}
and
\begin{align} \label{ubadatxk}
\lim_{i} u_i(P_{k,i}) = -\infty.
\end{align}

\end{corollary}

\begin{proof}
Taking $\theta = 1/2$ in Proposition \ref{lifenearx1} we get
\begin{align*}
u(x) \leq \ds \frac{3}{2}\log d_g(x,x_1) + C
\end{align*}
for $x$ near $x_1$.  On the other hand,(\ref{eregforu}) implies
\begin{align*}
u(x) \geq \log d_g(x,x_1) - C^{\prime}.
\end{align*}
Since these inequalities contradict one another when $x$ is close
enough to $x_1$, we conclude that (\ref{badcase}) is false.
Therefore, (\ref{supatpointbig}) must hold.

Now, according to Proposition \ref{nearsigv1}, given any compact $K
\subset M^n \setminus \Sigma$, there is a constant $C = C(K)>0$ such
that
\begin{align} \label{hrnk}
\max_{K} |\nabla u_i| \leq C(K)
\end{align}
for all $i \geq J = J(K)$.  Therefore, choosing a subsequence so
that $\lim_{i} u_i(x_0) = \limsup_{i} u_{i}(x_0) = +\infty$, the
gradient bound (\ref{hrnk}) implies that (\ref{infisbig}) holds for
any $N > 0$ and all $i$ sufficiently large, at least for $K \subset
M^n \setminus \Sigma$. Once again, however, by restricting to a
subsequence we may be reducing the singular set.

Near each point $x_k \in \Sigma = \{ x_1, \dots ,x_{\ell} \}$, there
are two possibilities to consider. First, suppose in a neighborhood
$V$ of $x_k$ we have
\begin{align} \label{nearxkgood}
u_i \geq -C.
\end{align}
Then the local $C^2$-estimate of Guan and Wang (Lemma
\ref{gradrest}) would imply
\begin{align*}
|\nabla^2 u_i| + |\nabla u_i|^2 \leq C^{\prime}
\end{align*}
in a neighborhood $V^{\prime} \subset V$.  Since $u_i \to +\infty$
pointwise on $M^n \setminus \Sigma$, (\ref{infisbig}) is valid on
any compact $K  \subset M^n \setminus \big(\Sigma \setminus \{ x_k
\} \big)$ for $i$ sufficiently large.  In this case, $x_k \notin
\Sigma_0$.

The alternative to (\ref{nearxkgood}) is that near $x_k \in \Sigma$
there is a sequence of points $\{ P_{k,i} \}$ satisfying
(\ref{badpoints}) and (\ref{ubadatxk}).  In this case, $x _k \in
\Sigma_0$.

Finally, note that the subsequence $\{ u_i \}$ can always be chosen
so that $\Sigma_0 \neq \emptyset$.  Otherwise, $\{ u_i \}$ would
have to be bounded from below near each point in $\Sigma$, and
consequently on all of $M^n$.  Combining the gradient estimate of
Guan-Wang with the fact that $-u_i \to +\infty$ pointwise on $M^n
\setminus \Sigma$, we would conclude that
\begin{align*}
\min_{M^n} u_i \to +\infty.
\end{align*}
But this contradicts the conclusion of Proposition \ref{floordrops}.
\end{proof}

\begin{remark}
Any sequence $\{ u_i \}$ of $\Gamma$-admissible solutions of
(\ref{hessF}) has a subsequence which satisfies the conclusions of
Corollary \ref{uitoplusinfty}, since the proof just relies on the
results of Section \ref{localestimates} and Lemma \ref{ricvolume}.
For the same reasons, Corollary \ref{uitoplusinfty} is valid for any
sequence $\{ u_i \}$ of solutions to $\tilde{\Psi}_{t_i}[u_i]= 0$.
\end{remark}

\section{The re-scaled sequence}
\label{rescaledsection}

Since $\{ u_i \}$ is diverging to $+\infty$ away from the singular
set $\Sigma_0$, we need to normalize the sequence if we hope to
extract a limit.  Let $x_0 \notin \Sigma_0$ again be a ``regular"
point, and define
\begin{align} \label{wdef}
w_i(x) = w_i(x) - w_i(x_0).
\end{align}
Using the properties of $\{ u_i \}$ derived in the preceding
section, we first show that (a subsequence of) $\{ w_i \}$ converges
off $\Sigma_0$ to a $C_{loc}^{1,1}$-limit.

\begin{proposition} \label{wifacts}
$(i)$  Given $r > 0$ small enough, let $M_r^n = M^n \setminus
\cup_{x_k \in \Sigma_0} B(x_k,r)$.  Then there is a constant $C >
0$, which is independent of $i$ and $r$, such that
\begin{align} \label{wiC2}
|\nabla^2 w_i|(x) + |\nabla w_i|^2(x) \leq Cr^{-2}
\end{align}
for $x \in M_r^n$ and all $i > J = J(r)$.

\vskip.1in  \noindent $(ii)$  The sequence $\{ w_i \}$ is bounded
above:
\begin{align} \label{wiabove}
\max_{M^n} w_i \leq C
\end{align}
for some constant $C > 0$.

\end{proposition}

\begin{proof}  $(i)$  As we saw in Corollary \ref{uitoplusinfty},
the original sequence $\{ u_i \}$ is diverging uniformly to
$+\infty$ on compact sets $K \subset M^n \setminus \Sigma_0$.  Let
$y \in M_{r}^{n}$; then the ball $B(y,r/2) \subset M^n \setminus
\Sigma_0$.  Therefore,
\begin{align} \label{limitonB}
\lim_{i} \inf_{B(y,r/2)} u_i = +\infty.
\end{align}
Applying the local estimate Lemma \ref{gradrest} on the ball
$B(y,r/2)$, we conclude that each $u_i$ satisfies
\begin{align*}
|\nabla^2 u_i|(x) + |\nabla u_i|^2(x) \leq C  \big( r^{-2} + e^{-2
\inf_{B(y,r/2)}u_i}\big)
\end{align*}
for all $x \in B(y,r/4)$.  Of course, since $u_i$ and $w_i$ only
differ by a constant, this implies
\begin{align}  \label{gradwionMr}
|\nabla^2 w_i|(x) + |\nabla w_i|^2(x) \leq C  \big( r^{-2} + e^{-2
\inf_{B(y,r/2)}u_i}\big)
\end{align}
for all $x \in B(y,r/4)$.  By (\ref{limitonB}), there is a $J =
J(y)$ such that
\begin{align*}
e^{-2 \inf_{B(y,r/2)}u_i} < r^{-2}
\end{align*}
for $i > J$.  Substituting this into (\ref{gradwionMr}) we get
\begin{align}  \label{boundwr}
|\nabla^2 w_i|(x) + |\nabla w_i|^2(x) \leq C r^{-2}
\end{align}
for all $x \in B(y,r/4)$ and $i > J = J(y)$.

The balls $\{ B(y,r/4) \}_{y \in M_{r}^{n}}$ define an open cover of
$M_{r}^{n}$, and since $M_{r}^{n}$ is compact we can extract a
finite subcover $M_r^n \subset \cup_{\nu = 1}^{N}
B(y_{\nu},r_{\nu}/4)$. Let $J = \max_{1 \leq \nu \leq N}J_{\nu}$.
For any $x \in M_r^n$, there is a ball with $x \in
B(y_{\nu},r_{\nu}/4)$, and inequality (\ref{boundwr}) is valid for
$i > J$.  This proves (\ref{wiC2}).

\vskip.2in \noindent $(ii)$  Since $w_i(x_0) = 0$, in view of the
bound (\ref{boundwr}) there must be a small ball $B(x_0,\rho_0)$ and
a constant $C > 0$ such that
\begin{align*}
\sup_{B(x_0,\rho_0)} w_i \leq C,
\end{align*}
for all $i \geq 1$.  This implies
\begin{align*}
\int_{B(x_0,\rho_0)} e^{\alpha w_i} dvol_g \leq C_{\alpha}
\end{align*}
for any $\alpha > 0$.  By Proposition \ref{lpmaxbound}, we obtain a
global bound:
\begin{align*}
\max_{M^n} w_i \leq C.
\end{align*}
\end{proof}

The next result summarizes the properties of the limit $w = \lim_i
w_i$.  Recall from the proof of Theorem \ref{intgrad} the definition
of the tensor $S_g = Ric - 2\delta \sigma_1(A_g)g$.  We let $S_w$
denote $S$ with respect to the limiting metric $g_w = e^{-2w}g$.

\begin{corollary} \label{wilim}
A subsequence of $\{ w_i \}$ converges on compact sets $K \subset
M^n \setminus \Sigma_0$ in $C^{1,\beta}(K)$, any $\beta \in (0,1)$.
Moreover, \vskip.1in

\noindent $(i)$ the limit $w = \lim_{i} w_i$ is in
$C_{loc}^{1,1}(M^n \setminus \Sigma_0)$. \vskip.1in

\noindent $(ii)$  The Hessian $\nabla^2 w(x)$ is defined at almost
every $x \in M^n$. \vskip.1in

\noindent $(iii)$ The tensor $S_w(x)$ is positive semi-definite at
almost every $x \in M^n$. \vskip.1in
\end{corollary}

\begin{proof}
Most of the statements are immediate consequences of Proposition
\ref{wifacts}, the Arzela-Ascoli Theorem, and the fact that
$w_i(x_0) = 0$. From Rademacher's Theorem, $\nabla^2 w$ is
well-defined almost everywhere (meaning the matrix of second
partials is well-defined almost everywhere), and $\nabla^2 w \in
L^{\infty}_{loc}$. Statement $(iii)$ follows from a standard
limiting argument using an integration by parts; we therefore omit
the details.
\end{proof}

\begin{lemma}
\label{wlimgrowth} We have the following estimates for $w$:
\begin{align}
\label{decay2nextnew} |\nabla w|_g (x) \leq \frac{C}{d_g(x,x_k)},
\end{align}
and for any $s>1$,
\begin{align}
\label{decay3nextnew} | \nabla w|_{C^{0,1}(A_g(r, s r))} \leq
\frac{C}{r^2},
\end{align}
where $A_g(r, sr)$ is the annulus in the metric $g$, and we take the
Lipschitz seminorm: that is,
\begin{align*}
\|f\|_{C^{0,1}(\Omega)} = \underset{x,y \in \Omega, x \neq y }
{\mbox{ sup }} \frac{| f(x) - f(x)|}{d_g(x,y)},
\end{align*}
for any domain $\Omega$.
\end{lemma}
\begin{proof}
The estimates (\ref{wiC2}) hold for the $w_i$, and since $w$ is the
$C^{1,1}$-limit obtained using the Arzela-Ascoli theorem, the lemma
follows immediately.
\end{proof}

\begin{remark}  The preceding analysis can be applied to any sequence of
$\Gamma$-admissible solutions $\{ u_i \}$ to (\ref{hessF}), rescaled
so as to converge in $C^{1,\alpha}$ on compact sets in the manner
described by Corollary \ref{wilim}.  In particular, the limit $w =
\lim w_i$ will satisfy the conclusions of Corollary \ref{wilim}.

In fact, all the results of the next two Sections \ref{singsection}
and \ref{finalsection} apply to such (rescaled) sequences of
solutions.  For this reason, from now on we will refrain from
calling this fact to the reader's attention.
\end{remark}

\section{Analysis of the singularities} \label{singsection}

The main result of this section is Theorem \ref{summagrowth}, which
gives a preliminary estimate of the limiting function near each $x_k
\in \Sigma_0$:
\begin{align} \label{wislike0}
2 \log d_g(x,x_k) - C \leq w(x) \leq 2 \log d_g(x,x_k) + C.
\end{align}
As we shall see, the first inequality above is a fairly easy
consequence of the maximum principle; the second inequality,
however, is much more delicate. One consequence of this estimate is
that the metric
\begin{align} \label{gcompletedef}
g_{w} = e^{-2w}g
\end{align}
in complete.  In Section 7 we will give further refinements of the
asymptotic behavior of $w$ near $\Sigma_0$, which in turn give us a
better understanding of the behavior of the metric $g_{w}$ near
infinity.

We begin with the proof of the first inequality in (\ref{wislike0}).
This estimate would be an easy consequence of the maximum principle,
except for the fact that $w$ is not $C^2$. Therefore, we need to
prove the corresponding statement for the $w_i$'s, then take a
limit.

\begin{proposition} \label{propwfromabove}
There is a constant $C$ such that
\begin{align} \label{wifromabove}
w_i(x) \geq 2 \log d_g(x,x_k) - C
\end{align}
for all $x$ near $x_k \in \Sigma_0$.
\end{proposition}

\begin{proof}

It will simplify matters if we use the notation introduced in the
proof of Proposition \ref{lpmaxbound}.  Let $v_i =
e^{-\frac{(n-2)}{2}w_i}$; then by (\ref{Scaleqnv}) $v_i$ satisfies
\begin{align} \label{Scaleqnvi}
Lv_i = \Delta v_i - \frac{(n-2)}{4(n-1)}Rv_i \leq 0,
\end{align}
where $L$ is the conformal laplacian.  Given $x_k \in \Sigma_0$, let
$\Gamma_k = \Gamma(x_k,\cdot)$ denote the Green's function for $L$
with singularity at $x_k$ (note that $L$ exists since $R > 0$).
Define
\begin{align} \label{Gammadef}
\Gamma = \sum_{k} \Gamma_k,
\end{align}
and
\begin{align} \label{Gdef}
G_i = \ds \frac{v_i}{\Gamma}.
\end{align}
Then an easy calculation using (\ref{Scaleqnvi}) shows
\begin{align} \label{Gpde}
\Delta G_i + 2\langle \nabla G_i, \nabla \Gamma / \Gamma \rangle
\leq 0
\end{align}
on $M^n \setminus \Sigma_0$, where $\langle \cdot, \cdot \rangle$
denotes the inner product. By the maximum principle, for all $r > 0$
small enough
\begin{align} \label{mpGi}
\min_{M_r^n} G_i = \min_{\cup \partial B(x_k,r) } G_i.
\end{align}
Now, by the gradient estimate (\ref{wiC2}),
\begin{align*}
\max_{\partial B(x_k,r)} w_i - \min_{\partial B(x_k,r)} w_i \leq C_n
r \max_{\partial B(x_k,r)} |\nabla w_i| \leq C_n r(C/r) = C,
\end{align*}
and consequently we have the Harnack inequality
\begin{align*}
\max_{\partial B(x_k,r)} v_i \leq C \min_{\partial B(x_k,r)} v_i,
\end{align*}
independent of $r$.  Since $\Gamma$ satisfies a similar inequality,
it follows that $G_i$ must:
\begin{align*}
\max_{\partial B(x_k,r)} G_i \leq C \min_{\partial B(x_k,r)} G_i.
\end{align*}
Therefore, by (\ref{mpGi}),
\begin{align} \label{supGi}
\max_{\cup \partial B(x_k,r)} G_i \leq C \min_{M_r^n} G_i.
\end{align}
The last inequality implies that $G_i$ is bounded on all of $M^n
\setminus \Sigma_0$.  To see this, first note that for $r > 0$
small, if $x_0$ is our regular point then (\ref{supGi}) implies
\begin{align*}
\max_{\cup \partial B(x_k,r)} G_i \leq C G_i(x_0) =
C\Gamma^{-1}(x_0) \leq C.
\end{align*}
This shows that
\begin{align*}
\max_{M^n \setminus \Sigma_0} G_i \leq C
\end{align*}
independent of $i$.  Therefore
\begin{align*}
\ds \frac{ e^{-\frac{(n-2)}{2}w_i} }{ \Gamma } \leq C.
\end{align*}
Finally, since
\begin{align*}
\Gamma(x_k,x) \sim d_g(x_k,x)^{(2-n)},
\end{align*}
by taking logarithms we get (\ref{wifromabove}).

\end{proof}

\begin{corollary} \label{corwfromabove}
There is a constant $C$ such that
\begin{align} \label{wfromabove}
w(x) \geq 2 \log d_g(x,x_k) - C
\end{align}
for all $x$ near $x_k \in \Sigma_0$.
\end{corollary}

\begin{proof} Inequality (\ref{wfromabove}) follows from
(\ref{wifromabove}) by letting $i \to \infty$ and using the fact
that $\{ w_i \}$ converges in $C^{1,\beta}$ off of $\Sigma_0$.
\end{proof}

To prove the second inequality in (\ref{wislike0}) we proceed in two
stages.  First, we prove a slightly weaker version of the result:

\begin{proposition}  \label{thetaclose}
Given $\theta > 0$, there is a constant $C = C(\theta)$ such that
near each $x_k \in \Sigma_0$,
\begin{align} \label{almostgood}
w(x) \leq (2 - \theta)\log d_g(x,x_k) + C(\theta)
\end{align}
for all $x \neq x_k$.
\end{proposition}

\begin{proof}  Using Corollary \ref{uitoplusinfty}, the proof of (\ref{almostgood}) is identical in its details
to the proof of Proposition \ref{lifenearx1} $(ii)$, and will
therefore be omitted.
\end{proof}

In the following, fix $x_k \in \Sigma_0$ and let $\rho = \rho(x) =
d_g(x,x_k)$.  We will assume that $\rho$ is well defined for $\rho <
\rho_0 < 1$.  An easy consequence of (\ref{almostgood}) is

\begin{lemma} \label{ratiois2}
\begin{align*}
\lim_{x \to x_k} \ds \frac{w(x)}{ \log \rho(x) } = 2.
\end{align*}
\end{lemma}

\begin{proof}
Given any $\theta > 0$, by (\ref{wfromabove}) and (\ref{almostgood})
we know
\begin{align*}
2 \log \rho(x) - C \leq w(x) \leq (2 - \theta)\log \rho(x) +
C(\theta)
\end{align*}
Therefore,
\begin{align*}
2 \geq \lim_{x \to x_k} \ds \frac{w(x)}{\log \rho(x)} \geq (2 -
\theta).
\end{align*}
Since $\theta > 0$ was arbitrary, Lemma \ref{ratiois2} follows.
\end{proof}

We now give the proof of the second inequality of (\ref{wislike0}):

\begin{theorem} \label{sharpgrowth}
Near $x_k \in \Sigma_0$, the function $w$ satisfies
\begin{align} \label{daineq}
w(x) \leq 2 \log \rho(x) + C.
\end{align}
\end{theorem}

\begin{proof}   As in the preceding proofs, the argument is
complicated by the fact that $w$ is not in $C^2$.   Moreover, we
cannot argue, as we did in the proof of Proposition
\ref{propwfromabove}, by establishing the result for the sequence
$\{ w_i \}$ and then passing to a limit: the $w_i$'s clearly do not
satisfy (\ref{daineq}).   To get around this difficulty we first
prove an estimate for $w$ on a dyadic annulus of fixed radius, then
iterate the estimate to bound the growth of $w$ near $x_k$.  The key
technical ingredient in this analysis is the following Proposition:

\begin{proposition} \label{thekeyprop}
Choose a point, say $x_1 \in \Sigma_0$, and let $B(a)$ denote the
geodesic ball of radius $a > 0$ centered at $x_1$. There exist
constants $K,a_0 > 0$ such that if $a < a_0$,
\begin{align} \label{keyestimate}
\max_{\partial B(a/2)} w \leq \max_{\partial B(a)} w - 2 \log a +
\Big[ 2 - K(a/2)^2\Big] \log (a/2).
\end{align}
\end{proposition}

\begin{proof}
For fixed $a > 0$ and small let
\begin{align} \label{bigWdef}
W(x) = w(x) - \max_{\partial B(a)}w + 2\log a,
\end{align}
and
\begin{align} \label{bigFdef}
F(x) = \ds \frac{ W(x) }{\log \rho}.
\end{align}
Later in the proof we will impose further conditions on $a$, but for
now we just require that $a < 1$ is small enough so that $\rho$ is
well defined in $B(a)$.  Note that $\log \rho(x) < 0$ for $x \in
B(a)$.

Now, by Lemma \ref{ratiois2} we have
\begin{align} \label{F(0)}
\lim_{x \to x_1} F(x) = 2.
\end{align}
Also, if $x \in \partial B(a)$, then
\begin{align*}
W(x) &= w(x) - \max_{\partial B(a)}w + 2\log a \\
&\leq 2 \log a,
\end{align*}
which implies
\begin{align} \label{F(a)}
\min_{\partial B(a)} F \geq 2.
\end{align}
From (\ref{F(0)}) and (\ref{F(a)}) we conclude that either
\begin{align} \label{Fbigger2}
F(x) \geq 2 \quad \forall x \in B(a),
\end{align}
or $F$ attains its minimum in the interior of $B(a) \setminus
\{x_1\}$.  However, if (\ref{Fbigger2}) holds then inequality
(\ref{keyestimate}) follows almost immediately.  To see this, let
$y_0$ be a point at which $F$ attains its minimum on $\partial
B(a/2)$.  Then (\ref{Fbigger2}) implies
\begin{align*}
2 \leq F(y_0) = \ds \frac{\max_{\partial B(a/2)} w -
\max_{\partial B(a)} w - 2\log a} {\log (a/2)} \\
\\
\Rightarrow \quad \max_{\partial B(a/2)} w \leq \max_{\partial B(a)}
w - 2 \log a + 2\log (a/2).
\end{align*}
Thus, (\ref{keyestimate}) holds with $K = 0$.  We may therefore
assume $F$ attains its minimum at a point $z_0 \in B(a) \setminus
\{x_1\}$.

To prove (\ref{keyestimate}) we want to apply the strong maximum
principle to $F$, though one needs to amend this statement slightly
because of regularity considerations.  As we observed above, $w \in
C_{loc}^{1,1}$, while the distance function $\rho$ is smooth on the
deleted ball $B(a) \setminus \{ x_1  \}$. Thus, $F$ is
$C_{loc}^{1,1}$ in a small neighborhood $U_0$ of $z_0$. Given
$\delta > 0$, we define a differential operator $\mathcal{L}$ in
$U_0$ by
\begin{align} \label{defL}
\mathcal{L}u = a^{ij}\nabla_i \nabla_j u,
\end{align}
where $\{a^{ij}\}$ are the components of the tensor
\begin{align} \label{defa}
a = (n-2)d\rho \otimes d\rho + (1 - 2\delta)g.
\end{align}
For $\delta > 0$ sufficiently small $\mathcal{L}$ is obviously
strictly elliptic. If $F$ were $C^2$ it would follow that
$\mathcal{L}F(z_0) \geq 0$, since $z_0$ is a minimum point. However,
$F$ is only $C_{loc}^{1,1}$, so while $\nabla^2 F$ exists almost
everywhere it may happen that $\mathcal{L}F(z_0)$ is not defined.
Despite this, there must be points nearby for which $\mathcal{L}F$
is {\it almost} non-negative:

\begin{lemma} \label{weakmaxF}
There is a sequence of points $\{z_j\}$ in $U_0$ satisfying
\vskip.2in \noindent $(i)$ $z_j \to z_0$ as $j \to \infty$,

\vskip.2in \noindent $(ii)$ $\nabla^2F(z_j)$ exists,

\vskip.2in \noindent $(iii)$ $\liminf_{j \to \infty}
\mathcal{L}F(z_j) \geq 0$,

\vskip.2in \noindent $(iv)$ $\nabla F(z_j) \to 0$ as $j \to \infty$.

\end{lemma}

\begin{proof}  Since $z_0$ is a minimum point of $F \in
C^{1,1}(U_0)$, conclusion $(iv)$ holds for {\it any} sequence of
points $z_j \to z_0$.  If $(iii)$ failed for every sequence
satisfying $(i)$ and $(ii)$, then there would be a neighborhood $U_1
\subseteq U_0$ such that
\begin{align} \label{supersoln}
\mathcal{L}F(z) \leq 0
\end{align}
at every point $z \in U_1$ where $\nabla^2 F(z)$ exists (in fact,
$\mathcal{L}F(z) < 0$ at every such point, but for our purposes it
is enough to assume the weaker inequality). In particular,
(\ref{supersoln}) would hold almost everywhere in $U_1$.

Now, $\mathcal{L}$ is not written in divergence form, but this is
easy to do:
\begin{align} \label{divformL}
\mathcal{L}u = \nabla_{i}\big(a^{ij}\nabla_j u \big) - b^j \nabla_j
u,
\end{align}
where $\{ b^j \}$ are the components of the divergence of
$\{a^{ij}\}$:
\begin{align*}
b^j = \nabla_k a^{jk}.
\end{align*}
It is important to note that the neighborhoods $U_1 \subseteq U_0$
do not contain $x_1$, so the distance function $\rho$, and
consequently the coefficients of $\mathcal{L}$, are all smooth on
$U_1$.

Returning to inequality (\ref{supersoln}), if $\varphi \in
C_{0}^{\infty}(U_1)$, then by (\ref{supersoln}) and (\ref{divformL})
we have
\begin{align*}
\int_{U_1} \Big( - a^{ij}\nabla_i \varphi \nabla_j F - \varphi b^j
\nabla_j F \Big) dvol \leq 0.
\end{align*}
Therefore, $F$ also satisfies (\ref{supersoln}) in $U_1$ in a weak
$W^{1,2}$-sense. By the strong maximum principle for weak
supersolutions (see \cite{GT}, Theorem 8.19), it follows that $F$
cannot have an interior minimum unless $F$ is constant on $U_1$.
Since $z_0$ is an interior minimum $F$ must be constant on $U_1$, in
which case the conclusions of Lemma \ref{weakmaxF} are obviously
true.

\end{proof}

We now apply Lemma \ref{weakmaxF} by calculating $\mathcal{L}F$ at
the points $\{ z_j \}$ described above.   To begin,
\begin{align} \label{gradF}
\nabla_l F = \ds \frac{ \nabla_l W }{ \log \rho} - \frac{ W}{(\log
\rho )^2} \frac{\nabla_l \rho }{ \rho}.
\end{align}
Differentiating again, and using (\ref{gradF}), we obtain
\begin{align} \label{hessF2}
\nabla_k \nabla_l F = \ds \frac{ \nabla_k \nabla_l W }{\log \rho} -
\frac {\nabla_k F }{\log \rho} \frac{\nabla_l \rho}{\rho} - \frac
{\nabla_l F }{\log \rho} \frac{\nabla_k \rho}{\rho}
  + \frac{F}{\log \rho}\frac{\nabla_k
\rho}{\rho}\frac{\nabla_l \rho}{\rho} - \frac{F}{ \log
\rho}\frac{\nabla_k \nabla_l \rho}{\rho}.
\end{align}
Tracing (\ref{hessF2}) we have
\begin{align} \label{lapF}
\Delta F = \ds \frac{ \Delta W}{ \log \rho} - \frac{2}{\log \rho}
\langle \nabla F, \frac{\nabla \rho}{\rho} \rangle + \frac{F}{\log
\rho} \frac{1}{\rho^2} - \frac{F}{\log \rho}\frac{ \Delta
\rho}{\rho}.
\end{align}

Using (\ref{hessF2}) and (\ref{lapF}) we can now write the
expression for $\mathcal{L}F$.  In doing so we will make use of the
following identities:
\begin{align} \label{altLdef}
\mathcal{L}u &= (n-2)\nabla^2 u(\nabla \rho,\nabla \rho) + (1 -
2\delta)\Delta u
\end{align}
for any function $u$, and
\begin{align*}
|\nabla \rho|^2 = 1,  \quad \nabla^2 \rho (\nabla \rho, \nabla \rho)
= 0.
\end{align*}
Thus,
\begin{align} \label{LFv1} \begin{split}
(\log \rho)\mathcal{L}F = \mathcal{L}W +& \ds
(n-1-2\delta)\frac{F}{\rho^2}
- (1-2\delta) F \frac{\Delta \rho}{\rho} \\
&-2(n-1-2\delta)\langle \nabla F, \frac{\nabla \rho}{\rho} \rangle.
\end{split}
\end{align}
Since the background metric $g$ is assumed to be $k$-admissible it
must have positive Ricci curvature.  Therefore, by the Laplacian
Comparison Theorem,
\begin{align*}
\ds \frac{\Delta \rho}{\rho} \leq \frac{(n-1)}{\rho^2}.
\end{align*}
Substituting this into (\ref{LFv1}) we get
\begin{align} \label{LFv2}
(\log \rho)\mathcal{L}F \geq \mathcal{L}W + \ds
2\delta(n-2)\frac{F}{\rho^2} - 2(n-1-2\delta)\langle \nabla F,
\frac{\nabla \rho}{\rho} \rangle.
\end{align}

The next step involves estimating $\mathcal{L}W$.  By Corollary
\ref{wilim}, for almost every $z$, the Ricci curvature satisfies
\begin{align} \label{Ricposw}
Ric(g_w) - 2\delta \sigma_1(A_w)g \geq 0
\end{align}
for $\delta \leq \frac{(2k - n)(n-1)}{2n(k-1)}$. In what follows, we
fix a value of $\delta > 0$ satisfying this condition. Rewriting
(\ref{Ricposw}) using (\ref{Sdef}) and (\ref{Schange}), we get
\begin{align} \label{WRic} \begin{split}
(n&-2)\nabla^2 W + (1-2\delta)\Delta W g + (n-2)dW \otimes dW \\
&- (n-2)(1 - \delta)|\nabla W|^2 g + \big(Ric_g -
2\delta\sigma_1(A_g)g\big) \geq 0. \end{split}
\end{align}
Using (\ref{gradF}), (\ref{altLdef}),
and the inequality $|\langle \nabla F, \nabla \rho \rangle| \leq
|\nabla F|$, we can rewrite this as
\begin{align} \label{LW2} \begin{split}
\mathcal{L}W &\geq \ds -\delta(n-2)\frac{F^2}{\rho^2} -2\delta(n-2)F
\log \rho \langle \nabla F, \frac{ \nabla \rho
}{\rho} \rangle \\
& \quad \quad - \delta(n-2)(1-\delta) (\log \rho)^2 |\nabla F|^2 -
C. \end{split}
\end{align}
Substituting this into (\ref{LFv2}), we finally arrive at
\begin{align*}
(\log \rho)\mathcal{L}F \geq \ds -\delta(n-2)\frac{F^2}{\rho^2} +
2\delta(n-2)\frac{F}{\rho^2} - C + \big(\mbox{terms with $\nabla
F$}\big).
\end{align*}
Keep in mind that we are evaluating both sides of the above
inequality at the points $\{z_j \}$ described in Lemma
\ref{weakmaxF}.  In particular, by property $(iv)$ of this lemma,
$|\nabla F(z_j)| \to 0$ as $j \to \infty$, and using property
$(iii)$ of the lemma along with the fact that $\rho(z_j) \to
\rho(z_0) = \rho_0$,
it follows that for all $j$ sufficiently large
\begin{align*}
C \geq \ds -\delta(n-2)\frac{F^2(z_j)}{\rho(z_j)^2} +
2\delta(n-2)\frac{F(z_j)}{\rho(z_j)^2}.
\end{align*}
Letting $j \to \infty$,  and completing the square, we obtain
\begin{align*}
\big(F(z_0) - 1\big)^2 \geq 1 - C\rho_0^2.
\end{align*}
Solving this inequality, since $\rho_0 \ll 1$, one sees that there
are two possible conclusions: either
\begin{align} \label{curtain1a}
F(z_0) \geq 2 - C_1 \rho_0^2
\end{align}
or
\begin{align} \label{curtain2a}
F(z_0) \leq C_2\rho_0^2 \leq C_2a^2.
\end{align}
We want to rule out (\ref{curtain2a}).  To do so we rely on
Corollary \ref{corwfromabove}, which implies
\begin{align*}
\max_{\partial B(a)} w \geq 2 \log a - C_3,
\end{align*}
and consequently
\begin{align*}
W(z_0) \leq w(z_0) + C_3.
\end{align*}
Therefore,
\begin{align} \label{F0big}
F(z_0) = \ds \frac{W(z_0)}{\log \rho_0} \geq \frac{w(z_0)}{ \log
\rho_0} - \frac{C_3}{\log \rho_0}.
\end{align}
From Lemma \ref{ratiois2},
we may then choose $a > 0$ small enough so that
\begin{align} \label{F1big}
F(z_0) > 1,
\end{align}
for all $x \in B(a)$. At the same time, by (\ref{curtain2a}), we can
choose $a > 0$ small enough (just depending on the constant $C_2$)
so that
\begin{align} \label{F2big}
F(z_0) \leq \frac{1}{2}.
\end{align}
Since (\ref{F1big}) and (\ref{F2big}) obviously contradict one
another, we conclude that (\ref{curtain1a}) must hold:
\begin{align*}
\min_{B(a)} F = F(z_0) \geq 2 - C_1 \rho_0^2 \geq 2 - C_1 a^2.
\end{align*}
Recalling the definition of $F$, we now have the following:
\begin{align*}
\ds \frac{ w(x) - \max_{\partial B(a)}w + 2\log a }{\log \rho(x)}
\geq 2 - C_1 a^2.
\end{align*}
If we take $x$ to be a point at which $w$ attains its maximum on
$\partial B(a/2)$, we obtain the inequality
\begin{align*}
\ds \frac{ \max_{\partial B(a/2)} w - \max_{\partial B(a)}w + 2\log
a }{\log (a/2)} \geq 2 - C_1 a^2 = 2 - K(a/2)^2.
\end{align*}
Multiplying both sides of the above inequality by $\log (a/2)$ and
rearranging terms, we arrive at (\ref{keyestimate}):
\begin{align} \label{keyagain}
\max_{\partial B(a/2)} w \leq \max_{\partial B(a)} w - 2 \log a +
\Big[ 2 - K(a/2)^2\Big] \log (a/2).
\end{align}
\end{proof}

To pass from Proposition \ref{thekeyprop} to Theorem
\ref{sharpgrowth}, we iterate inequality (\ref{keyagain}) to obtain
that for any $N > 1$,
\begin{align} \label{iter2} \begin{split}
\max_{\partial B(a/2^N) } w &\leq \max_{\partial B(a)} w - 2\log a
+ 2 \log (a/2^N) \\
&\quad \quad -K \Big[ (a/2)^2 \log (a/2) + \cdots + (a/2^N)^2 \log
(a/2^N) \Big].  \end{split}
\end{align}
Since the series
\begin{align*}
(a/2)^2 \log (a/2) + \cdots = \sum_{k=1}^{\infty} (a/2^k)^2 \log
(a/2^k)
\end{align*}
clearly converges, (\ref{iter2}) implies
\begin{align} \label{iter3}
\max_{\partial B(a/2^N) } w &\leq \max_{\partial B(a)} w - 2\log a +
2 \log (a/2^N) + C_4,
\end{align}
where $C_4 = C_4(a)$.

To complete the proof of the Theorem, let $x \in B(a)$, and choose
an integer $N \geq 1$ with $a/2^{N+1} < \rho(x) \leq a/2^N$.  By the
gradient estimate (\ref{decay2nextnew}),
\begin{align*}
w(x) - \max_{\partial B(a/2^N)} w &\leq
\Big(\frac{a}{2^N}\Big)\max_{B(a/2^N) \setminus B(a/2^{N+1})}
|\nabla w|
\leq C.
\end{align*}
Therefore, by (\ref{iter3}),
\begin{align*}
w(x) &\leq \max_{\partial B(a/2^N)} w + C \leq 2 \log(a/2^N) + C
\leq 2\log \rho(x) + C.
\end{align*}
This completes the proof.

\end{proof}

Combining Corollary \ref{corwfromabove} and Theorem
\ref{sharpgrowth}, we have

\begin{theorem} \label{summagrowth}
There is a constant $C$ such that
\begin{align} \label{wfromabovefin}
2 \log d_g(x,x_k) - C \leq w(x) \leq 2 \log d_g(x,x_k) + C
\end{align}
for all $x$ near $x_k \in \Sigma_0$.
\end{theorem}

\section{Volume growth estimates and completion of proof}
\label{finalsection}

In this section we refine our estimates of the asymptotic behavior
of $w$ near the singular set, with the eventual goal of showing that
the metric $g_w = e^{-2w}g$ is isometric to the Euclidean metric.

\subsection{Preliminary asymptotic behavior}

As before, we write $g_w = e^{-2w}g$, and the singular set is
denoted $\Sigma_0 = \{ x_1, \dots x_{\nu} \}$.  For $x_k \in
\Sigma_0$, let $r_{x_k} (x) = d_g(x_k, x)$ denote the distance in
the background metric. Note that $r_{x_k}^2$ is Lipschitz on $M^n$,
and is smooth in a neighborhood of $x_k$.

Let $\tilde{r}$ be a smooth function on $M^n$, positive on $M^n
\setminus \Sigma_0$, such that $\tilde{r} = r_{x_k}$ in a
neighborhood of each singular point $x_k \in \Sigma_0$. We write the
metric as
\begin{align}
g_w = e^{-2w} g = e^{-2 \Psi} \cdot \tilde{r}^{-4} g =  e^{-2 \Psi}
g_{ \star}
\end{align}
where
\begin{align} \label{Psidef}
\Psi =  w - 2 \log \tilde{r},
\end{align}
\begin{align} \label{gstardef}
g_{\star} = \tilde{r}^{-4} g.
\end{align}

In this notation, Theorem \ref{summagrowth} can be restated as

\begin{theorem}
\label{summagrowth2} There exists a constant $C_1 > 0$ so that
\begin{align}
 -C_1 \leq \Psi(x) \leq C_1,
\end{align}
for all $x \in M^n_{reg} = M^n \setminus \Sigma_0$.
\end{theorem}

Therefore, the metric $g_w$ is a bounded perturbation of the metric
$g_{\star}$. The asymptotic behavior of $g_{\star}$ is given by

\begin{lemma}[\cite{LeeandParker}]
\label{confnorm} For $x_k \in \Sigma_0$, let $\{ y^j \}$ be a
Riemannian exponential normal coordinate system for $g$ in a
neighborhood $U_k$ of $x_k$. Let $z^j = |y|^{-2} y^j$ denote
inverted normal coordinates on $U_k \setminus \{x_k \}$.  Then with
respect to these coordinates,

\vskip.1in \noindent $(i)$ $r_{x_k}(y) = |y| = \sqrt{(y^1)^2 + \dots
(y^n)^2}$.

\vskip.1in \noindent $(ii)$ Let $\rho(z) = |z| = \sqrt{(z^1)^2 +
\dots (z^n)^2}$, then $g_{\star} = \rho^4 g.$

\vskip.1in \noindent $(iii)$  $g_{\star}$ is asymptotically flat of
order two; i.e.,
\begin{align*}
(g_{\star})_{ij} = \delta_{ij} + O( \rho^{-2}), \
\frac{\partial}{\partial z^l}(g_{\star})_{ij} = O( \rho^{-3}), \
\frac{\partial^2}{\partial z^l \partial z^m} (g_{\star})_{ij} = O(
\rho^{-4}).
\end{align*}
\end{lemma}

\begin{remark}
the metric $g_{\star}$ is asymptotically flat of order $2$, with a
Euclidean end corresponding to each $x_k \in \Sigma_0$.
\end{remark}

The preceding Lemma implies the following expansion for the metric
$g_w$.
\begin{lemma}
\label{invertlem} Near each singular point, in inverted normal
coordinates $\{ z^i \}$ we have
\begin{align}
\label{decay1} e^{-2C_1} \delta_{ij} + O( \rho^{-2})
\leq & (g_w)_{ij} \leq  e^{2C_1} \delta_{ij} + O( \rho^{-2}),\\
\label{decay2} \frac{\partial}{\partial z^l}(g_w)_{ij} &= O(
\rho^{-1}).
\end{align}
Furthermore, for any $s>1$ and any $l = 1 \dots n$,
\begin{align}
\label{decay3} \| \frac{\partial}{\partial z^l}
(g_w)_{ij}\|_{C^{0,1}(A_0(r, s r))} = O( r^{-2}),
\end{align}
where $A_0(r,sr) = \{ z : r < \rho(z) < s r \}$, and $\| \cdot
\|_{C^{0,1}(A_0(r, s r))}$ denotes the Lipschitz seminorm on the
annulus with respect to the Euclidean metric.
\end{lemma}

\begin{proof} Inequalities (\ref{decay1})-(\ref{decay3}) follow from elementary
calculations using Lemma \ref{confnorm}, along with the estimates
for $w$ in contained in Lemma \ref{wlimgrowth} and Theorem
\ref{summagrowth2}.
\end{proof}

\subsection{Geodesic completeness}
Let us fix $x_0 \in M_{reg}^n$ and write
\begin{align*}
M^n \setminus \Sigma_0 = M^n_{reg} = M_0 \cup M_1,
\end{align*}
where $M_0$ is a compact smooth Riemannian manifold with boundary
containing $x_0$, and
\begin{align*}
M_1 \approx \amalg_1^{\nu} (\mathbf{R}^n \setminus B(0,R_0)),
\end{align*}
where $B(0,R_0)$ denotes the open ball of radius $R_0
> 0$ in $\mathbf{R}^n$.  On each component
we also have the inverted coordinates $\{ z^l \}$ as described in
Lemma \ref{confnorm}. We will call these components {\em{ends}},
denote them by $\{ N_j \}_{j=1}^{\nu}$, and implicitly use inverted
coordinates to identify each $N_j$ with $\mathbf{R}^n \setminus
B(0,R_0)$. In particular, the coordinates $\{ z^l \}$ on $N_j$
define a Euclidean metric $ds^2$ with Euclidean distance function
$|\cdot|_0$.  For example, given a point $p \in N_1$, in inverted
coordinates we write $p = (z^1,\dots,z^n)$, and $|p|_0 =
\mbox{dist}(0,p) = |z| = \rho(z)$.

 Since $g_w$ is $C^{1,1}$, the geodesic equation is a second order
system of ODEs with Lipschitz coefficients. The local existence and
uniqueness of solutions follows from a standard ODE theorem; see,
for example, \cite{Coddington} or \cite{Hartman}.  In particular,
the exponential map $exp_{x_0}: T_{x_0}M^n \rightarrow M^n$ is
defined in some neighborhood of the origin in $T_{x_0}M^n.$  The
next proposition says that, in fact, geodesics can be infinitely
extended.
\begin{proposition}
\label{comp1} The manifold $(M \setminus \Sigma_0, e^{-2w}g)$ is
geodesically complete. That is, for any $x_0 \in  M^n_{reg}$, the
exponential map $\exp_{x_0}: T_{x_0} M^n \rightarrow M^n_{reg}$ (of
the metric $g_w = e^{-2w} g$) is defined on all of $T_{x_0} M^n$.
\end{proposition}
\begin{proof}

Let $x_0 \in M^n_{reg}$, and let $\gamma(t)$ be a unit-speed
geodesic with $\gamma(0) = x_0$. Assume the maximal domain of
definition of $\gamma$ is $[0,T)$; we want to show that $T =
\infty$.

We will make use of the following property of geodesics: $\gamma:
[0,T) \rightarrow M^n_{reg}$ must leave every compact subset of
$M^n_{reg}$ as $t \rightarrow T$. That is, given any compact subset
$K \subset M^n_{reg}$, there exists a $t_K$ such that $\gamma(t) \in
M^n_{reg} \setminus K$ for all $t
> t_K$. For a simple proof of this fact, see \cite[page
109]{Petersen}.

Therefore, without loss of generality, we may assume there is a time
$0 < a < T$ such that $\gamma(t) \in M_1$ for $t \geq a$. Without
loss of generality, we may assume that $\gamma(t) \in N_1$ for $t
\geq a$. For $t \geq a$, by (\ref{decay1})
\begin{align}
\begin{split}
\label{estC2}
| \dot{\gamma}(t)|^2_{g_w} & = (g_w)_{ij}(\gamma(t)) \dot{\gamma}(t)^i \dot{\gamma}(t)^j \\
& \leq \Big( e^{2C_1} \delta_{ij} + O ( |\gamma(t)|_0^{-2}) \Big)  \dot{\gamma}(t)^i \dot{\gamma}(t)^j\\
& \leq  \Big( e^{2C_1} + \frac{C}{|\gamma(t)|_0^2} \Big) |
\dot{\gamma}(t)|^2_{0}.
\end{split}
\end{align}
Note that in (\ref{estC2}), and the inequalities which follow, we
are using the identification $M_1 = \mathbf{R}^n \setminus
B(0,R_0)$; hence for $t \geq a$ we have $\gamma(t),\dot{\gamma}(t)
\in \mathbf{R}^n$, and $|\cdot|_0$ denotes the Euclidean norm.
Similarly, using (\ref{decay1}) again we have
\begin{align*}
| \dot{\gamma}(t)|^2_{g_w} & \geq  \Big( e^{-2C_1} -
\frac{C}{|\gamma(t)|_0^2} \Big) | \dot{\gamma}(t)|^2_{0}.
\end{align*}

Now let $b \in (a,T)$, so $\gamma(b) \in M_1$.  Since $\gamma$ has
unit speed, the length of $\gamma([0,b])$ is given by
\begin{align} \label{L0b}
\begin{split}
b = L( \gamma([0,b]) = \int_0^b | \dot{\gamma}(t)|_{g_w} dt
&= \int_0^{a}  | \dot{\gamma}(t)|_{g_w} dt + \int_a^b  | \dot{\gamma}(t)|_{g_w} dt \\
&\geq a + \int_a^b  \Big( e^{-C_1} - \frac{C}{ | \gamma(t)|_0} \Big)
| \dot{\gamma}(t)|_0 dt.
\end{split}
\end{align}
By the definition of $a$,
\begin{align*}
|\gamma(a)|_0 = \min_{t \in [a,b]} |\gamma(t)|_0,
\end{align*}
so the integrand in (\ref{L0b}) can be estimated as
\begin{align*}
\Big( e^{-C_1} - \frac{C}{ | \gamma(t)|_0} \Big)  |
\dot{\gamma}(t)|_0 \geq \Big( e^{-C_1} - \frac{C}{|\gamma(a)|_0}
\Big) | \dot{\gamma}(t)|_0.
\end{align*}
Substituting this into (\ref{L0b}) we obtain
\begin{align}  \label{dcomp1}
b \geq a + \Big( e^{-C_1} - \frac{C}{|\gamma(a)|_0} \Big) \int_a^b |
\dot{\gamma}(t)|_0 dt.
\end{align}
Since segments minimize distance in the Euclidean metric, we have
\begin{align*}
\int_a^b | \dot{\gamma}(t)|_0 dt \geq |\gamma(b) - \gamma(a)|_0.
\end{align*}
Therefore,
\begin{align}
\label{goody} b - a & \geq \Big( e^{-C_1}  - \frac{C}{|\gamma(a)|_0}
\Big) | \gamma(b) - \gamma(a)|_0.
\end{align}

Now, recall that given any compact set $K \subset M^n_{reg}$, there
must be a time $t_K$ with $\gamma(t) \in M^n_{reg} \setminus K$ for
$t
> t_K$. Therefore, by choosing a large enough compact set we can
arrange so that $\gamma(b) \in M^n_{reg} \setminus K$ and
$|\gamma(b) - \gamma(a)|_0$ is as large as we like.  By
(\ref{goody}), this means we can choose $b$ as large as we like,
i.e., $T = \infty$. It follows that $(M^n_{reg},g_w)$ is
geodesically complete.
\end{proof}

\subsection{Properties of the distance function}

 The distance function $d_w$ is defined on
$M^n_{reg} \times M^n_{reg}$ by
\begin{align}
d_w (p,q) = \inf_{\gamma} \{ L_w (\gamma), \gamma: [t_0,t_1]
\rightarrow M^n_{reg} \mbox{ is a piecwise } C^1 \mbox{ path joining
} p \mbox{ and } q \},
\end{align}
and
\begin{align}
L_w(\gamma) = \int_{t_0}^{t_1} | \dot{\gamma}(t)|_{g_w} dt.
\end{align}

We begin with a preliminary lemma on the convergence of the distance
functions.

\begin{lemma}
\label{distconv} For any compact subset $K \subset M^n_{reg}$,
\begin{align*}
\lim_{i \rightarrow \infty} d_{g_i} = d_w.
\end{align*}
uniformly on $K \times K$.
\end{lemma}
\begin{proof}
Take $p, q \in K$. First we show that
\begin{align} \label{deq1}
\limsup_i d_i(p,q) \leq d_w(p,q).
\end{align}
Choose $\epsilon > 0$, and let $\gamma_w$ denote a piecewise $C^{1}$
curve from $p$ to $q$ in $M^n_{reg}$ with $L_w(\gamma_w) \leq
d_w(p,q) + \epsilon$. If $L_{g_i}[\gamma_w]$ denotes the length of
$\gamma_w$ with respect to $g_i$, then
\begin{align} \label{deq2}
d_{i}(p,q) \leq L_{g_i}[\gamma_w].
\end{align}
Letting $i \to \infty$, since $g_i \to g$ on compact sets,
\begin{align}
\limsup_i d_i(p,q) \leq L_w[\gamma_w] \leq d_w(p,q) + \epsilon.
\end{align}
Since $\epsilon$ is arbitrary, we obtain (\ref{deq1}).

We now turn to the (more difficult) opposite inequality.  To this
end, let $\gamma_i: [\alpha_i,\beta_i] \rightarrow M^n_{reg}$ denote
a minimizing geodesic from $p$ to $q$ in the metric $g_i$.

\begin{claim} For $i >> 1$ sufficiently large, there is a compact set $K
\subset M^n_{reg}$ containing $p$ and $q$, such that the image of
$\gamma_i$ lies entirely in $K$.
\end{claim}
\begin{proof}
To verify this claim, we argue by contradiction.  Write $\Sigma = \{
x_1,\dots,x_{\nu} \}$, and let
\begin{align}
M_r = M^n \setminus \big( \cup_{k=1}^{\nu} B_g(x_k,r) \big),
\end{align}
where $B_g$ means a ball in the background metric $g$.  Suppose
that, for every $r > 0$, there is an infinite number of points $x_i
= \gamma_i(t_i)$ with $x_i \in \partial M_{r}$.  Since $\{w_i\}$
converges in $C^{1,\alpha}$ on $M_r$, inequality (\ref{compdiw})
holds on $M_r$ for every tangent vector $X$ and all sufficiently
large $i$.  Therefore,
\begin{align} \label{dbig} \begin{split}
d_i(p,x_i) =
\int_{\alpha_i}^{t_i}\sqrt{g_i(\dot{\gamma_i},\dot{\gamma_i})} dt
 &\geq (1 - \delta)^{1/2}
\int_{\alpha_i}^{t_i}\sqrt{g_w(\dot{\gamma_i},\dot{\gamma_i})}
dt \\
&\geq (1 - \delta)^{1/2} d_w(p,x_i) \\
&\geq (1 - \delta)^{1/2}d_w(p,\partial M_r), \end{split}
\end{align}
for $i$ sufficiently large, where $d_w(p,\partial M_r)$ denotes the
distance from $p$ to $\partial M_r$ in the metric $g_w$.

From equation (\ref{goody}), it follows that
\begin{align*}
d_w(p,\partial M_r) \rightarrow \infty
 \ \mbox{ as } r \to 0.
\end{align*}
In particular, (\ref{dbig}) tells us that by first choosing $r > 0$
small enough, we can make $d_i(p,x_i)$ as large as we like for $i$
sufficiently large.  More precisely:  Given $D
>> 1$, there is a $J = J(D)$ such that
$d_i(p,x_i) \geq D$ for all $i > J$.  But this implies
\begin{align*}
d_i(p,q) \geq d_i(p,x_i) \geq D,
\end{align*}
for $i > J$.  If $D$ is chosen large enough (say, $D = 2d_w(p,q)$),
this obviously contradicts (\ref{deq1}).

\end{proof}

Now to complete the proof, since $w_i \rightarrow w$ in
$C^{1,\alpha}(K)$, given $\delta > 0$ there is a $J >> 1$ such that
for all tangent vectors $X$,
\begin{align}  \label{compdiw} \begin{split}
g_i(X,X) &= e^{-2w_i}g(X,X) \\
&=  e^{-2(w_i - w)}g_w(X,X) \geq (1 - \delta)g_w(X,X)
\end{split}
\end{align}
for all $i > J$. Therefore,
\begin{align*}
d_i(p,q) = L_{g_i}[\gamma_i] &= \int_{\alpha_i}^{\beta_i}
\sqrt{g_i(\dot{\gamma_i},\dot{\gamma_i})} dt \\
&\geq (1 - \delta)^{1/2} L_w[\gamma_i] \geq (1 - \delta)^{1/2}
d_w(p,q).
\end{align*}
Taking the limit we obtain
\begin{align*}
\liminf_i d_i(p,q) \geq (1-\delta)^{1/2} d_w(p,q).
\end{align*}
Since this holds for all $\delta > 0$, we have
\begin{align*}
\liminf_i d_i(p,q) \geq d_w(p,q).
\end{align*}
\end{proof}

\begin{proposition}[Bishop volume comparison]
\label{growthless} Fix any basepoint $x_0 \in M^n_{reg}$. Then the
ratio
\begin{align}
\label{bishop}
 \frac{ Vol_{g_w}(B_{g_w}(x_0,r))}{ r^n} \leq \omega_n,
\end{align}
and is a non-increasing function of $r$, where $\omega_n$ is the
Euclidean volume growth constant.
\end{proposition}
\begin{proof}

Let $ g_{j} = e^{-2w_j} g$, where $\{w_j\}$ is the sequence
constructed in Section \ref{rescaledsection}.  Recall by Corollary
\ref{wilim} the metrics $\{ g_{j} \}$ converge to $g_w = e^{-2w}g$
in the $C^{1,\alpha}$-norm, for any $\alpha < 1$, on compact subsets
of $M^n_{reg}$. Let $d_j(x) = d_{w_j}(x_0, x)$ denote the distance
from $x_0$ in the metric $g_j$, and let $d (x) = d_w(x_0,x)$.

 Since each $g_{j}$ is a smooth metric with positive Ricci
curvature, inequality (\ref{bishop}) holds for each $g_j$ (see
\cite[Lemma 9.1.6]{Petersen}). That is, if $r_2 \geq r_1$ then
\begin{align} \label{Volj}
Vol_{g_j}(B_{g_j}(x_0,r_2)) r_2^{-n} \leq
Vol_{g_j}(B_{g_j}(x_0,r_1)) r_1^{-n} \leq \omega_n.
\end{align}
\begin{claim}The volume of balls satisfies
\begin{align} \label{volconv}
\lim_{j \rightarrow \infty} Vol_{g_j}(B_{g_j}(x_0,r)) =
Vol_{g_w}(B_{g_w}(x_0,r))
\end{align}
as $j \to \infty$.
\end{claim}
\begin{proof}
This follows directly from Lebesgue's dominated convergence theorem.
\end{proof}

To complete the proof of the Proposition we take the limit in
(\ref{Volj}).
\end{proof}

\begin{lemma}
\label{Cheeger} For each point $x_0 \in M^n_{reg}$, there exists a
radius $r_{x_0}$ so that the distance function $d_w(x_0, \cdot) \in
C^{1,1}( B_{g_w}( x_0, r_{x_0}) \setminus \{x_0\})$. Furthermore,
the exponential map in the $g_w$ metric, $exp_{x_0}: T_{x_0} M \cap
B(0, r_{x_0}) \rightarrow B_{g_w}(x_0, r_{x_0})$, is a Lipschitz
homeomorphism.
\end{lemma}
\begin{proof}
By Proposition \ref{wifacts} and Corollary \ref{wilim}, the sequence
$\{ w_j \}$ and its derivatives up to order two are uniformly
bounded on $\overline{B_0} = \overline{B_w(x_0,r_0)}$, for some
radius $r_0$. It follows that the metrics $\{ g_j \}$ have uniformly
bounded curvature on $\overline{B_0}$:
\begin{align} \label{curvbound}
\max_{\overline{B_0}} |Riem_{g_j}| \leq C.
\end{align}

If we take $r_2 = r_0$ and $r_1 = r$ in (\ref{Volj}), then
\begin{align*}
Vol_{g_j}(B_{g_j}(x_0,r_0)) r_0^{-n} \leq Vol_{g_j}(B_{g_j}(x_0,r))
r^{-n}.
\end{align*}
Letting $j \to \infty$ and using (\ref{volconv}) implies
\begin{align*}
Vol_{g_j}(B_{g_j}(x_0, r_0)) \rightarrow Vol_{g_w}(B_{g_w}(x_0,
r_0)) > 0.
\end{align*}
 In particular, we have a lower bound
\begin{align} \label{lowvol}
Vol_{g_j}(B_{g_j}(x_0,r)) \geq v_0 \cdot r^n
\end{align}
for some $v_0 > 0$ and all $r < r_0$.

From the curvature bound (\ref{curvbound}) and the volume bound
(\ref{lowvol}), by \cite[Theorem 4.7]{CGT} we conclude that the
injectivity radius $inj_{g_j}(x_0)$ of $g_j$ at $x_0$ is strictly
bounded from below, with the estimate
\begin{align*}
inj_{g_j}(x_0) \geq \Lambda r_0> 0,
\end{align*}
for some constant $\Lambda > 0$. In particular, for $r < \Lambda
r_0$, the distance function $d_j^2 = d_j(x_0, \cdot)^2$ is smooth
inside of $B_{g_j}(x_0, r)$.

  Furthermore, using the curvature
bound (\ref{curvbound}), from \cite[Lemma 10.4.2]{Petersen}, we
obtain a bound $ | \nabla^2_j d_j^2| \leq C$ on a possibly smaller
ball $B_{g_j}(x_0, \Lambda' r_0)$. We have shown in Lemma
\ref{distconv} that $d_j \rightarrow d_w$ on compact subsets. From
Arzela-Ascoli, by passing to another subsequence if necessary, we
have that and $d_j^2 \rightarrow d_w^2$ in $C^{1, \alpha}$, and that
$d_w^2 \in C^{1,1}$ in $B_{g_w}(x_0, \Lambda' r_0)$. This proves the
first statement.

  For the statement about the exponential map, the Lipschitz
differentiability follows from the standard ODE dependence on
initial conditions, see \cite[Chapter V, Theorem 8.1]{Hartman}. The
exponential map must be injective on $B(0, r_{x_0})$; this follows
since $d^2(x_0, \cdot)$ is $C^{1,1}$ there. To see why, supppose
$v_1, v_2 \in B(0, r_{x_0}) \subset T_{x_0} M$ satisfied
$\exp_{x_0}(v_1) = \exp_{x_0}(v_2) = q$. Then there would be $2$
distinct geodesics starting from $x_0$ ending at $q$. Without loss
of generality, we may assume $q$ is the first intersection point
after $x_0$. From local uniqueness of geodesics through $q$
mentioned above, we would then have two distinct radial directions
at $q$, which would imply that the distance function is not even
$C^1$ at $q$, a contradiction.
\end{proof}

 \begin{proposition}[Hopf-Rinow]
\label{Hopf-Rinow} For each $x \in M^n_{reg}$, the exponential map
in the $g_w$ metric, $\exp_x: T_x M^n \rightarrow M^n_{reg}$ is
onto. Therefore, any $x_1, x_2 \in M^n_{reg}$ can be joined by a
$C^2$-geodesic $\gamma : [0, d_w(x_1,x_2)] \rightarrow M^n_{reg}$.
\end{proposition}
\begin{proof}
The local unique minimizing property of geodesics follows from Lemma
\ref{Cheeger} (see \cite[Lemma 5.3.6]{Petersen}), and from the
identity $|\nabla_w d_w| = 1$ on $B(x, r_x) \setminus \{x\}$.
However, the local minimizing property, together with the properties
of the exponential map proved in Lemma \ref{Cheeger}, are the key
ingredients in the standard proof of Hopf-Rinow (see, for example,
\cite[Theorem 5.7.1]{Petersen}).
\end{proof}

\subsection{Distance function estimates}
 In the previous section, we showed that the manifold
$(M^n_{reg},g_w)$ is complete, but we need more precise information
about the behavior of $g_w$ at infinity. We next show that for
points in the same end of $M^n_{reg}$, distances measured with
respect to $g_w$ are comparable to distances measured in the induced
Euclidean metric.

Recall that $M^n = M_0 \cup M_1$, with
\begin{align*}
M_1 \approx \amalg_1^{\nu} (\mathbf{R}^n \setminus B(0,R_0)).
\end{align*}
We denote the ends by $\{ N_j \}_{j=1}^{\nu}$, and as before use
inverted coordinates to identify each $N_j$ with $\mathbf{R}^n
\setminus B(0,R_0)$.

\begin{proposition} \label{comp3}
Fix $x_0 \in M_0$ as in Proposition \ref{comp1} and an end, say
$N_1$. Then there are constants $R_1 > 0$ and $C_2 > 0$ such that
for any $x \in N_1$ with $d_w(x_0,x) \geq R_1$, we have the distance
estimates
\begin{align}
\label{distcomp} e^{-C_2} |x|_0 - C_2 \leq d_w(x_0, x) \leq e^{C_2}
|x|_0 + C_2.
\end{align}
\end{proposition}

\begin{proof}
 We will establish (\ref{distcomp}) through a
series of technical lemmas. We begin with an estimate which was
proved in Proposition \ref{comp1}.

\begin{lemma} \label{compcor}
Let $p_1, p_2$ be in the same component of $M_1$ and assume that the
geodesic between $p_1$ and $p_2$ in the $g_w$ metric lies entirely
in $M_1$. Let $\gamma(a) = p_1$, and $\gamma(b) = p_2$, and assume
that
\begin{align}
\label{show} |\gamma(a)|_0 = \min_{t \in [a,b]} |\gamma(t)|_0.
\end{align} Then
\begin{align}
\Big( e^{-C_1}  - \frac{C}{|p_1|_0} \Big) | p_2 - p_1|_0 \leq
d_w(p_1, p_2).
\end{align}
\end{lemma}
\begin{proof}
This is a consequence of the formula (\ref{goody}).
\end{proof}

\begin{lemma}
\label{distlem1} There exist constants $R_1 > 0$ and $C_3 > 0$ such
that for $x \in N_1$ with $d_w(x_0,x) \geq R_1$, we have the
distance estimate
\begin{align}
\label{distest1}  e^{-C_3} |x|_0 - C_3 \leq d_w(x_0, x).
\end{align}
\end{lemma}
\begin{proof}

Let $x \in N_1$.  From Proposition \ref{Hopf-Rinow}, there exists a
unit-speed minimizing geodesic $\gamma$ with $\gamma(0) = x_0$ and
$\gamma( d_w(x_0,x)) = x$.  Choose $R^{\prime} < |x|_0$ so that
$B_w(x_0,R^{\prime}) \supset M_0$, where $B_w(x_0,R^{\prime})$
denotes the geodesic ball in the $g_w$ metric (the existence of such
an $R'$ also follows from Proposition \ref{Hopf-Rinow}).   Let
\begin{align} \label{t0def}
t_0 = \max \{ t : |\gamma(t)|_0 = R^{\prime}, \gamma(t) \in N_1 \
\mbox{for} \ t
> t_0 \}.
\end{align}
Let $p_0 = \gamma(t_0)$; then by Lemma \ref{compcor}
\begin{align}
\Big( e^{-C_1}  - \frac{C}{R^{\prime}} \Big) | x - p_0|_0 \leq
d_w(p_0, x).
\end{align}
If $|x|_0$ is sufficiently large, say $|x|_0 > R_1 \gg 1$, then
$R^{\prime}$ can be chosen large enough so that
\begin{align*}
e^{-C_1} - \frac{C}{R^{\prime}} \geq e^{-2C_1}.
\end{align*}
Thus,
\begin{align} \label{dc1}
e^{-2C_1}|x - p_0|_0 \leq d_w(p_0,x).
\end{align}

From the triangle inequality,
\begin{align*}
d_w( p_0, x) \leq d_w(p_0, x_0) + d_w(x_0, x) \leq R' + d_w(x_0, x).
\end{align*}
Substituting this into (\ref{dc1}),
\begin{align*}
d_w(x_0, x) &\geq d_w(p_0,x) - R^{\prime} \\
&\geq e^{-2C_1} | x - p_0|_0 - R' \\
&\geq e^{-2C_1}(|x|_0  - R^{\prime}) - R',
\end{align*}
which implies
\begin{align}
e^{-C_3}|x|_0  - C_3 \leq d_w(x_0, x)
\end{align}
for some $C_3 = C_3(C_1,R^{\prime})$.
\end{proof}

\begin{lemma}
Let $p_2 \in N_1$ with $|p_2|_0 \geq R_1$, where $R_1$ is given in
Lemma \ref{distlem1}. Let $p_1 = \frac{R_1}{|p_2|_0} p_2 \in N_1$.
Then
\begin{align}
d_w(p_1, p_2) \leq \Big( e^{C_1} + \frac{C}{ R_1} \Big) | p_2 -
p_1|_0.
\end{align}
\end{lemma}
\begin{proof}  Let $\gamma$ denote the line segment in $N_1$ joining $p_1$ and
$p_2$, and assume $\gamma$ has unit speed (as measured in the
Euclidean metric).   Thus, $\gamma: [0,|p_2 - p_1|_0] \rightarrow
N_1$ with $|\dot{\gamma}(t)|_0 = 1, \gamma(0) = p_1$, and
$\gamma(|p_2 - p_1|_0) = p_2.$  Then the length of $\gamma$ in the
$g_w$ metric is
\begin{align*}
L_{g_w}(\gamma) = \int_{0}^{|p_2 - p_1|_0}
|\dot{\gamma}(t)|_{g_w}dt.
\end{align*}
Clearly,
\begin{align*}
|p_1|_0 = |\gamma(0)|_0 = \min_{0 \leq t \leq |p_2 - p_1|_0}
|\gamma(t)|_0.
\end{align*}
By the estimate (\ref{estC2}),
\begin{align*}
L_{g_w}( \gamma) &\leq \int_0^{|p_2 - p_1|_0} \Big( e^{C_1} +
\frac{C}{ | \gamma(t)|_0} \Big) | \dot{\gamma}(t)|_{0} dt
\\
& \leq \Big( e^{C_1} + \frac{C}{ | p_1|_0} \Big)
\int_0^{|p_2 - p_1|_0} | \dot{\gamma}(t)|_{0} dt \\
& =  \Big( e^{C_1} + \frac{C}{R_1} \Big) |p_2 - p_1|_0.
\end{align*}
Since the distance is the infimum of the length of all paths, we
obtain
\begin{align*}
d_w(p_1, p_2) \leq \Big( e^{C_1} + \frac{C}{R_1} \Big) |p_2 -
p_1|_0.
\end{align*}

\end{proof}
%

\begin{lemma}
\label{distlem2} There is a constant $C_4 > 0$ such that for any $x
\in N_1$, we have the distance estimate
\begin{align}
\label{distest2} d_w(x_0, x) \leq e^{C_4} |x|_0 + C_4.
\end{align}
\end{lemma}%

\begin{proof}
Let $x \in N_1$, then $|x|_0 \geq R_1$.  Let $p_1 =
\frac{R_1}{|x|_0} x$, and apply the previous Lemma (with $p_2 = x$)
to get
\begin{align*}
d_w(p_1, x) & \leq \Big( e^{C_1} + \frac{C}{ R_1 } \Big) | x - p_1|_0\\
& =  \Big( e^{C_1} + \frac{C}{R_1} \Big) \Big(1 - \frac{R_1}{|x|_0}
\Big)
 | x|_0\\
& \leq  \Big( e^{C_1} + \frac{C}{ R_1 } \Big) | x|_0.
\end{align*}
Applying the triangle inequality
\begin{align*}
d_w(x_0, x) \leq d_w(x_0, p_1) + d_w (p_1,x),
\end{align*}
we obtain
\begin{align*}
d_w(x_0,x) \leq  d_w(x_0, p_1) + \Big( e^{C_1} + \frac{C}{ R_1 }
\Big) | x|_0.
\end{align*}

By construction, $|p_1|_0 = R_1$. From Proposition \ref{Hopf-Rinow},
there is a geodesic ball of radius $R_2$ such that $p_1 \in B_w(x_0,
R_2)$ (in particular, the radius of this ball does {\it not} depend
on the point $x$).   Therefore,
\begin{align*}
d_w(x_0,x) \leq  R_2 + \Big( e^{C_1} + \frac{C}{ R_1 } \Big) | x|_0.
\end{align*}
Choosing $C_4 = C_4(C_1,R_1,R_2)$ large enough, this implies
(\ref{distest2}).
\end{proof}
To complete the proof of Proposition \ref{comp3}, we combine
inequalities (\ref{distest1}) and (\ref{distest2}).  By choosing
$C_2 = C_2(C_3,C_4)$ large enough, inequality (\ref{distcomp})
follows.

\end{proof}
\subsection{Tangent cone analysis}

Propositions \ref{comp1} and \ref{comp3} give preliminary
information about the geometry of the ends $N_j$ by looking at the
behavior of the distance function near infinity.  We now turn to
results which describe the volume growth of $g_w$, which in turn
will give us information about the curvature decay.

The first result in this direction is an integral estimate for the
conformal factor $\Psi$ defined in (\ref{Psidef}).  Recall that $g_w
= e^{-2\Psi}g_{\star}$, and that $g_{\star}$ is asymptotically flat
of order two (by Lemma \ref{confnorm}). To pass from information
about $g_{\star}$ to information about $g_w$ requires estimating the
conformal factor $\Psi$.  This is the motivation of the following
result:

\begin{theorem}
\label{gdecay}  For each end $N_j, 1 \leq j \leq \nu$, we have
\begin{align}
\label{Lncond} \int_{N_j} | \nabla_{g_{\star}} \Psi |_{g_{\star}}^n
dvol_{g_{\star}} < \infty.
\end{align}
\end{theorem}

\begin{proof}  Given radii $0 < r_1 < r_2$, let
$A^{\star}(r_1,r_2)$ denote the annulus centered at $x_0 \in
M^n_{reg}$ with inner radius $r_1$ and outer radius $r_2$.  Recall
that $g_w = e^{-2\Psi}g_{\star}$, and from Theorem
\ref{summagrowth2} we know that $|\Psi| \leq C_1$.

Applying Theorem \ref{intgrad} with background metric $g_{\star}$,
$\delta = 0,$ and $p = n$,
\begin{align}
\label{localintgrad7}
\begin{split}
\int_{A^{\star}(r_1,r_2)} |\nabla_{g_{\star}} \Psi |^n e^{\alpha
\Psi} & dvol_{g_{\star}} \leq C\Bigg(
\int_{A^{\star}(\frac{1}{2}r_1,2r_2)} |Ric_{g_{\star}}|^{n/2}
e^{\alpha \Psi}
dvol_{g_{\star}} \\
&+ r_1^{-n} \int_{A^{\star}(\frac{1}{2}r_1,r_1)}e^{\alpha \Psi}
dvol_{g_{\star}} + r_2^{-n} \int_{A^{\star}(r_2, 2r_2)} e^{\alpha
\Psi} dvol_{g_{\star}} \Bigg),
\end{split}
\end{align}
for any fixed $\alpha > 0$.  We claim that the integrals on the
$RHS$ of (\ref{localintgrad7}) remain bounded with $r_1 > 0$ fixed
and as $r_2 \to \infty$.  Once we verify this, (\ref{Lncond})
follows.

By Lemma \ref{confnorm}, the metric $g_{\star}$ is asymptotically
flat of order two.  This implies that the curvature and volume
growth of $g_{\star}$ satisfy
\begin{align} \label{Ricdecay}
\max_{A^{\star}(r,2r)} |Ric_{g_{\star}}| \leq Cr^{-4},
\end{align}
\begin{align} \label{Volstar}
Vol_{g_{\star}}\big(A^{\star}(r,2r)\big) \leq Cr^n,
\end{align}
for some constant $C$ and all large $r \gg 1$.  Therefore,
\begin{align*}
\int_{A^{\star}(\frac{1}{2}r_1,2r_2)} |Ric_{g_{\star}}|^{n/2}
e^{\alpha \Psi} dvol_{g_{\star}} \leq Cr_1^{-n},
\end{align*}
\begin{align*}
r_2^{-n} \int_{A^{\star}(r_2, 2r_2)} e^{\alpha \Psi}
dvol_{g_{\star}} \leq C.
\end{align*}
Thus, (\ref{Lncond}) holds.
\end{proof}

\begin{corollary} \label{Lnw}
For each end $N_j, 1 \leq j \leq \nu$, we have
\begin{align}
\label{Lnwe} \int_{N_j} | \nabla_{g_w} \Psi |_{g_w}^n dvol_{g_w} <
\infty.
\end{align}
\end{corollary}

\begin{proof}  This is a consequence of the conformal invariance
of the $L^n$-norm of the gradient.  That is, since $g_w =
e^{-2\Psi}g_{\star}$,
\begin{align*}
| \nabla_{g_w} \Psi |_{g_w}^n dvol_{g_w} =  | \nabla_{g_{\star}}
\Psi |_{g_{\star}}^n dvol_{g_{\star}}.
\end{align*}
\end{proof}

Using Theorem \ref{gdecay} and Corollary \ref{Lnw}, we can now study
the volume growth of $g_w$.

\begin{proposition}
\label{annulargrowth} Fix an integer $m > 1$ and an end, say $N_1$.
Let $x_0 \in M_0$ be the same point as in the statements of
Propositions \ref{comp1} and \ref{comp3}. Then
\begin{align} \label{asmvol}
\lim_{r \rightarrow \infty} Vol_{g_w}( A_1(m^{-1}r, r)) r^{-n}  =
\omega_n ( 1 - m^{-n}),
\end{align}
where $A_1(m^{-1}r,r) = \{ x \in N_1 | m^{-1}r \leq d_w(x_0, x) \leq
r \}$ is the annulus in the metric $g_w = e^{-2w}g$ in the end
$N_1$, and $\omega_n$ is the volume ratio of Euclidean space.
\end{proposition}
\begin{proof}
The proof is based on a construction of the geometric tangent cone
corresponding to the end $N_1$, as follows. Take any sequence of
numbers $r_i$ such that $r_i \rightarrow \infty$ as $i \rightarrow
\infty$, and consider the sequence of annuli $A_1(2^{-1} r_i, 2 r_i)
= \{ x \in N_1 | 2^{-1}r_i \leq d_w(x_0, x) \leq 2r_i \}$. From now
on we will drop the subscript, and it will be understood that for
all annuli under consideration we take the component that lives in
the end $N_1$.

For $r_i$ large enough there are inverted coordinates $\{ z^j \}$
defined on
\begin{align*}
A_i = A(2^{-1} r_i,2 r_i) \subset N_1.
\end{align*}
We introduce the rescaled coordinates $(z^j)^{\prime} =
\frac{z^j}{r_i}$ and the rescaled metrics $\tilde{g}_i = r_i^{-2}
g_w$, both of which are defined on the rescaled annuli
\begin{align*}
\big( \tilde{A_i}(2^{-1}, 2), \tilde{g}_i \big) = (A_i,
r_i^{-2}g_w).
\end{align*}
The inequality (\ref{distcomp}) implies that $\tilde{A_i}(2^{-1},
2)$ are bounded domains in $\mathbf{R}^n$, and are contained in a
Euclidean annulus $A_0 = A_0(\epsilon,\epsilon^{-1})$ of bounded
size. Using the estimates (\ref{decay1}), (\ref{decay2}), and
(\ref{decay3}), by the Arzela-Ascoli theorem it follows that for
some subsequence (which we continue to index by $i$), $\tilde{g}_i
\rightarrow g_{\infty}$ on compact subsets of $A_0$ in $C^{1,
\alpha}$, for any $\alpha < 1$.

Now, by the scale-invariance of the estimate in Corollary \ref{Lnw},
\begin{align*}
\lim_{i \to \infty} \int_{\tilde{A_i}} | \nabla_{\tilde{g}_i} \Psi
|_{\tilde{g}_i}^n dvol_{\tilde{g}_i} =  \lim_{i \to \infty}
\int_{A_i} | \nabla_{g_w} \Psi |_{g_w}^n dvol_{g_w} = 0.
\end{align*}
It follows that $g_{\infty}$ must be a constant times the Euclidean
metric.  In particular, the limiting metric is flat.

\begin{claim}
The distance function to the basepoint $d_w = d_w(x_0, \cdot)$ is a
viscosity solution of the Hamilton--Jacobi equation $| \nabla_w d_w
|_w = 1$, away from the base point $x_0$.
\end{claim}
\begin{proof}
First we observe that each $d_{w_i}= d_{w_i}(x_0, \cdot)$ is a
viscosity solution of the Hamilton-Jacobi equation $|\nabla_i
d_{w_i}|_{w_i} = 1$ away from $x_0$ (see \cite{Mantegazza},
\cite{Azagra}). Using the convergence of the $d_i$ from Lemma
\ref{distconv}, and convergence of the metrics in $C^{1,\alpha}$,
the claim then follows from the well-known fact that a uniform limit
of viscosity solutions to a sequence of uniformly convergent
equations is a viscosity solution of the limiting equation (see
\cite[Theorem 1.4]{CEL}, or \cite[Theorem I.2]{CrandallLionsfirst}).
Alternatively, using  Lemma \ref{Cheeger} and the Hopf-Rinow
Proposition \ref{Hopf-Rinow}, the method of proof of \cite[Theorem
6.23]{Azagra} shows directly that $d_w$ is a viscosity solution of
$| \nabla_w d_w |_w = 1$.
\end{proof}

 Next, consider the rescaled distance functions $\rho_i(z')
= r_i^{-1} d_w(p_0, r_i z')$.  These are Lipschitz continuous
functions on $\tilde{A_i}$ which also satisfy $|\nabla_{g_i} \rho_i|
= 1$ in the viscosity sense. Again, by the Arezla-Ascoli theorem,
some subsequence converges on compact subsets of $A_0$ in
$C^{\alpha}$, for any $\alpha < 1$, to a limit $\tilde{\rho}$. Since
the subsequence $\rho_i$ converges to a limit, the domains
$\tilde{A_i} = \tilde{A_i}(2^{-1}, 2)$ converge to a limiting domain
$\tilde{A}( 2^{-1}, 2)$. Note that $\tilde{\rho}$ is Lipschitz
continuous, and is a viscosity solution of $|\nabla_{g_{\infty}}
\tilde{\rho}| = 1$ on the domain $\tilde{A}( 2^{-1}, 2)$.

\begin{claim}
For $z' \in \tilde{A}( 2^{-1}, 2)$ we have
\begin{align}
\label{distcompr}
 e^{-C_2} | z'| \leq \tilde{\rho} (z') \leq e^{C_2} |z'|
\end{align}
\end{claim}
\begin{proof}
Given $z' \in  \tilde{A}( 2^{-1}, 2)$, let $z_i = r_i z'$. The
distance estimates from Proposition \ref{comp3} are
\begin{align*}
e^{-C_2} |z_i|_0 - C_2 \leq d_w(x_0, z_i) \leq e^{C_2} |z_i|_0 +
C_2.
\end{align*}
Dividing by $r_i$, we obtain
\begin{align*}
e^{-C_2} |z'|_0 - C_2 r_i^{-1} \leq r_i^{-1} d_w(x_0, z_i) \leq
e^{C_2} |z'|_0 + C_2 r_i^{-1}.
\end{align*}
As $i \rightarrow \infty$,
\begin{align*}
e^{-C_2} |z'|_0 \leq \lim_{i \rightarrow \infty} r_i^{-1} d_w(x_0,
z_i) \leq e^{C_2} |z'|_0.
\end{align*}
Since
\begin{align*}
 \lim_{i \rightarrow \infty} r_i^{-1} d_w(x_0, r_i z') =
\lim_{i \rightarrow \infty} \rho_i ( z') = \tilde{ \rho}(z'),
\end{align*}
the claim follows.
\end{proof}

For the subsequence of radii $r_i$ chosen above we repeat this
procedure for each integer $l >2$, with annuli of size $A( l^{-1}
r_i, l r_i)$, and obtain a nested sequence of limiting domains
\begin{align*}
\tilde{A}( 2^{-1}, 2) \subset \tilde{A}( 3^{-1}, 3) \subset \dots
\subset \tilde{A}( l^{-1}, l) \subset \cdots.
\end{align*}
Clearly, the estimate (\ref{distcompr}) holds on every $\tilde{A}(
l^{-1}, l)$, from which it easily follows that
\begin{align*}
\bigcup_{l=2}^{\infty} \tilde{A}( l^{-1}, l) = \mathbf{R}^n
\setminus \{0 \}
\end{align*}
with the flat metric.  Furthermore, the rescaled distance functions
$\rho_i$ converge on compact subsets of $\mathbf{R}^n \setminus {0}$
to $\tilde{\rho}$ in $C^{\alpha}$ for any $\alpha < 1$.  The
inequalities (\ref{distcompr}) imply that any level set $ \{
\tilde{\rho} = t \}$ is contained in a Euclidean annulus $A_{Euc}(
e^{-2C_2}t,  e^{C_2} t) \subset \mathbf{R}^n \setminus \{0 \}$. This
implies that by defining $\tilde{\rho}(0) = 0$, $\tilde{\rho}$ has a
continuous extension to $\mathbf{R}^n$, and $\tilde{\rho}$ is a
Lipschitz viscosity solution of the equation $|\nabla_0
\tilde{\rho}| = 1$, on $ \mathbf{R}^n \setminus \{0 \}$. By a
standard uniqueness result for solutions of Hamilton-Jacobi
equations (\cite[Theorem IV.2.6]{Bardi}; see also \cite[Theorem
3.1]{Mantegazza}, \cite[Theorem 6.23]{Azagra}), we must have
$\tilde{\rho}(x) = |x|_0$, the Euclidean distance to the origin.

To complete the proof of the Proposition, given any sequence $r_i
\rightarrow \infty$, the annuli $A(m^{-1}r_i, r_i) \subset A( m^{-1}
r_i, m r_i)$, and the rescalings of the larger annulus converge to
Euclidean annuli.  It necessarily follows that the rescalings of the
smaller annuli also converge to Euclidean annuli, and therefore
(\ref{asmvol}) holds.
\end{proof}

\begin{remark}
Note that we have proved more than claimed--not only does the volume
ratio converge, but all rescaled annuli converge to Euclidean
annuli. In geometric terms, we have shown there is a unique tangent
cone at infinity, which is Euclidean space.
\end{remark}

\begin{corollary} \label{growthgreater} For each end $N_j$,
\begin{align} \label{supereuclid}
\lim_{r \rightarrow \infty} Vol_{g_w}( B(x_0, r) \cap N_j)r^{-n}
\geq \omega_n.
\end{align}
\end{corollary}
\begin{proof}
Fix an integer $m > 0$. We have $A \left( m^{-1} r, r \right)
\subset B(x_0,r)$, so
\begin{align}
Vol_{g_w} \left( A \left( m^{-1} r, r \right) \right) r^{-n} \leq
Vol_{g_w} ( B(x_0,r) ) r^{-n}.
\end{align}
Taking the limit, and using Proposition \ref{annulargrowth},
\begin{align}
 \omega_n ( 1 - m^{-n}) \leq
\lim_{r \rightarrow \infty} Vol_{g_w} ( B(x_0,r) ) r^{-n}.
\end{align}
This holds for any $m > 0$, so the corollary follows by letting $m
\rightarrow \infty$.
\end{proof}

\subsection{Completion of proof}

By Corollary \ref{growthgreater}, the asymptotic growth of the
volume of geodesic balls in the metric $g_w$ is Euclidean.  If $g_w$
were a (smooth) metric of non-negative Ricci curvature, the equality
case of the Bishop comparison theorem would imply that $g_w$ is
flat. The goal of this section is to show that a ``weak" version of
Bishop's theorem remains valid, and therefore we can still conclude
$g_w$ is flat.  We do not attempt to proof a general result for
$C^{1,1}$ metrics, but instead rely on the special fact that $g_w$
is the limit of smooth metrics with positive Ricci curvature.

At the end of the section we complete the proof of Theorem
\ref{compsoln} by showing that, if $g_w$ is the Euclidean metric,
then $(M^n,g)$ is conformally equivalent to the round sphere.

\begin{proposition}
\label{isoexist} The function $w \in C^{\infty}(M^n_{reg})$, is
smooth.  That is, the metric $g_w = e^{-2w} g$ is smooth and there
exists a smooth isometry
\begin{align*}
\Phi: (M^n_{reg}, g_w) \rightarrow ( \mathbf{R}^n, g_{Euc}),
\end{align*}
where $g_{Euc}$ denotes the Euclidean metric.
\end{proposition}
\begin{proof}
Fix any point $x_0 \in M_{reg}$. Combining Corollary
\ref{growthgreater} and Proposition \ref{growthless}, we see that
there can be at most one end; that is, $M^n_{reg} = M^n \setminus
\{x_1\}$, and furthermore that
\begin{align}
\label{Voleq} Vol(B_{g_w}(x_0,r)) = \omega_n \cdot r^n
\end{align}
for all $r>0$.

Recall that we have shown in Lemma \ref{Cheeger} that $d_w^2 \in
C^{1,1}$ in $B_{g_w}(x_0, r_{x_0})$, and $d_j^2 \rightarrow d_w^2$
in $C^{1, \alpha}$. Inside the cutlocus, $d_j^2$ is smooth, so by
the Laplacian comparison theorem we have
\begin{align}
\Delta_j d_j^2 \leq 2n.
\end{align}
Since $d_j$ converges to $d_w$, this implies that the Laplacian
comparison theorem
\begin{align}
\label{lct} \Delta_w d_w^2 \leq 2n
\end{align}
holds for $g_w$ in the weak $W^{1,2}$-sense. Since $d_w^2 \in
C^{1,1}(B_{g_w} ( x_0, r_{x_0}))$, Rademacher's Theorem implies that
\begin{align}
\Delta_w  d_w^2 \leq 2n \ \mbox{a.e. in }  B_{g_w}(x_0, r_{x_0}).
\end{align}

Next, fix $0 < r_1 < r_2 < r_{x_0}$. Integrating by parts, we have
\begin{align}
\int_{A_{g_w}(r_1,r_2)} \Delta_w d_w^2 dV_w = 2 \{ r_2
\mathcal{H}^{n-1}( S( r_2)) - r_1 \mathcal{H}^{n-1} (S(r_1)) \},
\end{align}
Where $\mathcal{H}^{n-1}(S(r))$ denotes the surface area (with
respect to the metric $g_w$) of the submanifold $S(r) = \{ x : d_w
(x_0,x) = r \}$. Since $d_w$ in $C^{1,1}$ inside of
$A_{g_w}(r_1,r_2)$, this formula follows from the divergence
theorem, and is also valid for $C^{1,1}$ metrics. We then obtain
\begin{align*}
2n \omega_n ( r_2^n - r_1^n) &= 2n Vol( A_{g_w}(r_1, r_2))
\geq \int_{A_{g_w}(r_1,r_2)} \Delta_w d_w^2 dV_w \\
& =  2 \{ r_2 \mathcal{H}^{n-1}( S( r_2)) - r_1 \mathcal{H}^{n-1}
(S(r_1)) \}.
\end{align*}
Since $d_w$ is $C^{1,1}$, differentiating (\ref{Voleq}) implies that
\begin{align*}
\mathcal{H}^{n-1}( S( r)) = n \omega_n r^{n-1}.
\end{align*}
Substituting in the above, we find that
\begin{align*}
\int_{A_{g_w}(r_1,r_2)} \Delta_w d_w^2 dV_w = 2n \omega_n ( r_2^n -
r_1^n).
\end{align*}

We conclude that $d_w$ is a $C^{1,1}$ solution of the equation
\begin{align}
\label{hum1} \Delta_{g_w} d_w^2 = 2n \ \ \mbox{a.e. in }
B_{g_w}(x_0, r_{x_0}) .
\end{align}

 We next write out the above equation with respect to the
background metric. Recall the following formulas for the conformal
change of the Hessian:
\begin{proposition}[\cite{Jeff3}]
\label{confprop} For any function h,
\begin{align}
\label{hessianformula} {\nabla}^2_w h= \nabla^2_g h + dw \otimes dh
+ dh \otimes dw - \langle dw, dh \rangle_g g,
\end{align}
\begin{align}
\label{Lapformula}
 \Delta_w h = e^{-2w} \left( \Delta_g h + (2-n)  \langle dw, dh \rangle_g
\right).
\end{align}
\end{proposition}
Equation (\ref{hum1}) therefore is
\begin{align}
\label{hum1conf} \Delta_{g}d_w^2 + (2-n) \langle dw, d(d_w)^2
\rangle_g = 2ne^{2w} \ \ \mbox{a.e. in }  B_{g_w}(x_0, r_{x_0}) .
\end{align}
Notice that the second term on the left hand side is in $W^{1,p}$
for any $p < \infty$. By a standard regularity theorem (see
\cite{GT}), we conclude that $d_w^2 \in W^{3,p} (B_{g_w}(x_0,
r_{x_0})$ for any $1 < p < \infty$.

Since $ |\nabla d_w|_{g_w} = 1$, equation (\ref{hum1}) is equivalent
to
\begin{align}
\label{hum2} \Delta_{g_w} d_w = \frac{n-1}{d_w}.
\end{align}

\begin{claim} The distance function $d_w$ satisfies
\begin{align}
\label{kolp} \langle \nabla \Delta d_w, \nabla d_w \rangle +
|\nabla^2 d_w|_{g_w}^2 = - Ric_{g_w}(\nabla d_w, \nabla d_w),
\end{align}
a.e. in $B_{g_w}(x_0, r_{x_0}) \setminus \{x_0 \}$.
\end{claim}
\begin{proof}
We recall the Bochner-Lichnerowicz formula (\cite{Lichnerowicz}): If
$F \in C^{3}$ is any function, and $h \in C^2$ is any Riemannian
metric, then
\begin{align}
\label{lich} \frac{1}{2} \Delta_h | \nabla F|_h^2 = |\nabla^2 F
|_h^2 + \langle \nabla \Delta_h F, \nabla F \rangle_h + Ric_h(
\nabla F, \nabla F).
\end{align}
By an approximation argument, forumla (\ref{lich}) is valid for any
$F \in W^{3,p}$, $p \gg 1$. If $h$ is of the form $h = e^{-2 v} g$,
for a smooth function $v$, then using Proposition \ref{confprop}, we
see that (\ref{lich}) involves at most second derivatives of $v$. By
another approximation argument, (\ref{lich}) holds for $v \in
C^{1,1} = W^{2, \infty}$. In particular, we can take $F = d_w$ on
$B_{g_w}(x_0, r_{x_0}) \setminus \{x_0 \}$, and $h = g_w = e^{-2w}
g$, and we are done.
\end{proof}

Using the inequality
\begin{align}
\label{matrixineq} \frac{1}{n-1} ( \Delta d_w)^2 \leq |\nabla^2
d_w|^2,
\end{align}
and equation (\ref{hum2}), we obtain
\begin{align}
\begin{split}
\label{squeeze} 0 &= \langle \nabla \left( \frac{n-1}{d_w} \right),
\nabla d_w \rangle
+ \frac{1}{n-1} (\Delta d_w)^2 \\
&\leq  \langle \nabla \Delta d_w, \nabla d_w \rangle + |\nabla^2
d_w|^2 =  - Ric_{g_w}(\nabla d_w, \nabla d_w).
\end{split}
\end{align}
From Corollary \ref{wilim}, the Ricci curvature of $g_w$ is
non-negative almost everywhere. Therefore (\ref{squeeze}) implies
that
\begin{align}
 Ric_{g_w} (\nabla d_w, \nabla d_w) = 0,
\end{align}
almost everywhere in $B_{g_w}(x_0, r_{x_0})$.

 From the inequality $Ric_{g_w} \geq 2 \delta \sigma_1 (A_{g_w})g_w$
of Corollary \ref{wilim}, we find that $R_{g_w} = 0$ almost
everywhere in $B_{g_w}(x_0, r_{x_0})$. Consequently $v =
e^{-\frac{(n-2)}{2}w}$ satisfies
\begin{align}
\Delta_g v - \frac{(n-2)}{4(n-1)}R_g v = 0
\end{align}
almost everywhere in $B_{g_w}(x_0, r_{x_0})$. From elliptic
regularity (see \cite[Chapter 9]{GT}), it follows that $w \in
C^{\infty}(B_{g_w}(x_0, r_{x_0}) )$.

  We can repeat the above argument for any other
basepoint $x_0$, to conclude the metric $g_w$ is smooth in
$M_{reg}$. By the standard version of Bishop's theorem for smooth
metrics, it follows that $g_w$ is isometric to the Euclidean metric
(see \cite[Theorem 3.9]{Chavel}).
\end{proof}

Theorem \ref{compsoln} follows almost immediately. To see this,
Proposition \ref{isoexist} implies, in particular, that $(M^n_{reg},
g)$ is locally conformally flat, so from continuity, $(M^n, g)$ is
locally conformally flat (by vanishing of the Cotton tensor in
$n=3$, or vanishing of the Weyl tensor in $n
> 3$; see \cite{Eisenhart}). We also see that $M^n \setminus
\{ x_1 \}$ is diffeomorphic to $\mathbf{R}^n$, so $M^n \setminus \{
x_1 \}$ is the one point compactification of $\mathbf{R}^n$ and
$M^n$ is homeomorphic to $S^n$. In particular, $M^n$ is simply
connected and locally conformally flat, so by Kupier's Theorem
\cite{Kuiper}, $(M^n,g)$ is conformally equivalent to $(S^n,
g_{round})$, a contradiction.

 Note that the same arguments give the proof of Theorem
\ref{compsoln2}, based on the Remark at the end of Section
\ref{rescaledsection}.  Also, by the results of Section
\ref{degreeF}, we have completed the proof of Theorem \ref{Main2}.

\section{Appendix} \label{App}

\subsection{On the $A^{\tau}$-problem}

In this Section we make a few remarks about Example $3$ from the
introduction.
\begin{theorem}
Let $A^{\tau} \in \Gamma_k^+$, $k > n/2$, and
\begin{align}
{\tau} > {\tau}_0(n,k) = \frac{2(n-k)}{n}.
\end{align}
Then there exists a constant $\delta_0(n,k,{\tau}) > 0$ so that
\begin{align}
Ric > 2\delta_0(n,k,{\tau}) \sigma_1(A) \cdot g.
\end{align}
\end{theorem}
\begin{proof}
It follows from \cite{GVW} that if $A^{\tau} \in \Gamma_k^+$, then
\begin{align}
Ric \geq  \Bigg( {\tau} + \frac{(k-n)(2n - n{\tau} -2)}{n(k-1)}
\Bigg) \frac{R}{2(n-1)} g,
\end{align}
and the theorem follows by simple calculation.
\end{proof}
It follows that the operator $F$ in equation (\ref{sigmakt})
satisfies the desired properties.

\subsection{Uniform ellipticity of (\ref{hessF})}

Suppose $u \in C^2$ is a solution of (\ref{hessF}) with
\begin{align} \label{C2A}
\| u \|_{L^{\infty}} +  \| \nabla u \|_{L^{\infty}} + \| \nabla^2 u
\|_{L^{\infty}} \leq C.
\end{align}
We claim that (\ref{hessF}) is uniformly elliptic. In particular, by
the results of Evans \cite{Evans} and Krylov \cite{Krylov}, $u \in
C^{2,\alpha}$ and the Schauder estimates give classical regularity.

To prove the claim amounts to verifying the following property of
$F$:

\vskip.15in $(iii)^{\prime}$ If $\lambda \in \Gamma$ satisfies $(a)$
$F(\lambda) \geq c_1 > 0$, and $(b)$ $|\lambda| \leq L$, then there
is a constant $\epsilon_0 = \epsilon_0(c_1,L)$ such that
\begin{align} \label{UEA}
\epsilon_0^{-1} \geq \frac{\partial F}{\partial \lambda_i}(\lambda)
\geq \epsilon_0.
\end{align}
\vskip.15in

Note that if $u \in C^2$ is a solution satisfying the bound
(\ref{C2A}), then the eigenvalues of $A_u$ will obviously satisfy
properties $(a)$ and $(b)$. Moreover, inequality (\ref{UEA}) says
that the linearized operator is uniformly elliptic with ellipticity
constant $\epsilon_0$.

To prove property $(iii)^{\prime}$, let $\hat{\Gamma} = \{ \xi \in
\overline{\Gamma} : |\xi| = 1\}$, and $\hat{\lambda} =
\lambda/|\lambda| \in \hat{\Gamma}$.  We first show there is a
constant $\delta_0 > 0$ such that
\begin{align} \label{posdistA}
\mbox{dist}(\hat{\lambda},\partial \hat{\Gamma}) \geq \delta_0.
\end{align}
Assuming this for the moment, let us see how (\ref{UEA}) follows.

Since $F$ is homogenous of degree one, $\partial F/\partial
\lambda_i$ is homogenous of degree zero.  Therefore,
\begin{align*}
\frac{\partial F}{\partial \lambda_i}(\lambda) = \frac{\partial
F}{\partial \lambda_i}(\hat{\lambda}).
\end{align*}
But (\ref{posdistA}) says that $\hat{\lambda}$ is a fixed distance
from the boundary of $\hat{\Gamma}$.  Since $\partial F/\partial
\lambda_i > 0$ in $\Gamma$, inequality (\ref{UEA}) follows from the
continuity of $\partial F/\partial \lambda_i$.

Returning to (\ref{posdistA}), we argue by contradiction.  Suppose
there is a sequence $\{ \lambda_i \} \subset \Gamma$ with
$|\lambda_i| \leq L$ and $F(\lambda_i) \geq c_1 > 0$, but with
\begin{align} \label{posdistA2}
\mbox{dist}(\hat{\lambda_i},\partial \hat{\Gamma}) \rightarrow 0
\end{align}
as $i \rightarrow \infty$, where $\hat{\lambda_i} =
\lambda_i/|\lambda_i|$. Choose a subsequence (still denoted $\{
\hat{\lambda_i} \}$) with $\hat{\lambda_i} \rightarrow
\hat{\lambda}_0 \in \partial \hat{\Gamma}$.  By continuity,
$F(\hat{\lambda_i}) \rightarrow 0$. By homogeneity,
\begin{align}  \label{dont}
F(\lambda_i) = F(|\lambda_i|\hat{\lambda_i}) =
|\lambda_i|F(\hat{\lambda_i}) \rightarrow 0,
\end{align}
since $|\lambda_i| \leq L$.  However, (\ref{dont}) contradicts
assumption $(a)$.

\bibliography{Ricci_references}

\providecommand{\bysame}{\leavevmode\hbox to3em{\hrulefill}\thinspace}
\providecommand{\MR}{\relax\ifhmode\unskip\space\fi MR }
\providecommand{\MRhref}[2]{%
  \href{http://www.ams.org/mathscinet-getitem?mr=#1}{#2}
}
\providecommand{\href}[2]{#2}
\begin{thebibliography}{AFLM05}

\bibitem[AFLM05]{Azagra}
Daniel Azagra, Juan Ferrera, and Fernando Lopez-Mesas, \emph{{Nonsmooth
  analysis and Hamilton-Jacobi equations on Riemannian manifolds}}, J. Funct.
  Anal. \textbf{220} (2005), no.~2, 304--361.

\bibitem[BC64]{Bishop}
Richard~L. Bishop and Richard~J. Crittenden, \emph{Geometry of manifolds},
  Academic Press, New York, 1964.

\bibitem[BCD97]{Bardi}
Martino Bardi and Italo Capuzzo-Dolcetta, \emph{Optimal control and viscosity
  solutions of {H}amilton-{J}acobi-{B}ellman equations}, Systems \& Control:
  Foundations \& Applications, Birkh\"auser Boston Inc., Boston, MA, 1997, With
  appendices by Maurizio Falcone and Pierpaolo Soravia.

\bibitem[Bes87]{Besse}
Arthur~L. Besse, \emph{Einstein manifolds}, Springer-Verlag, Berlin, 1987.

\bibitem[CEL84]{CEL}
M.~G. Crandall, L.~C. Evans, and P.-L. Lions, \emph{Some properties of
  viscosity solutions of {H}amilton-{J}acobi equations}, Trans. Amer. Math.
  Soc. \textbf{282} (1984), no.~2, 487--502.

\bibitem[CGT82]{CGT}
Jeff Cheeger, Mikhail Gromov, and Michael Taylor, \emph{Finite propagation
  speed, kernel estimates for functions of the {L}aplace operator, and the
  geometry of complete {R}iemannian manifolds}, J. Differential Geom.
  \textbf{17} (1982), no.~1, 15--53.

\bibitem[CGY02a]{CGY2}
Sun-Yung~A. Chang, Matthew~J. Gursky, and Paul Yang, \emph{An a priori estimate
  for a fully nonlinear equation on four-manifolds}, J. Anal. Math. \textbf{87}
  (2002), 151--186, Dedicated to the memory of Thomas H.\ Wolff.

\bibitem[CGY02b]{CGY1}
Sun-Yung~A. Chang, Matthew~J. Gursky, and Paul~C. Yang, \emph{An equation of
  {M}onge-{A}mp\`ere type in conformal geometry, and four-manifolds of positive
  {R}icci curvature}, Ann. of Math. (2) \textbf{155} (2002), no.~3, 709--787.

\bibitem[Cha93]{Chavel}
Isaac Chavel, \emph{Riemannian geometry---a modern introduction}, Cambridge
  Tracts in Mathematics, vol. 108, Cambridge University Press, Cambridge, 1993.

\bibitem[Che05]{Sophie1}
Szu-yu~Sophie Chen, \emph{Local estimates for some fully nonlinear elliptic
  equations}, preprint, 2005.

\bibitem[CL55]{Coddington}
Earl~A. Coddington and Norman Levinson, \emph{Theory of ordinary differential
  equations}, McGraw-Hill Book Company, Inc., New York-Toronto-London, 1955.

\bibitem[CL83]{CrandallLionsfirst}
Michael~G. Crandall and Pierre-Louis Lions, \emph{Viscosity solutions of
  {H}amilton-{J}acobi equations}, Trans. Amer. Math. Soc. \textbf{277} (1983),
  no.~1, 1--42.

\bibitem[CNS85]{CNSIII}
L.~Caffarelli, L.~Nirenberg, and J.~Spruck, \emph{The {D}irichlet problem for
  nonlinear second-order elliptic equations. {I}{I}{I}. {F}unctions of the
  eigenvalues of the {H}essian}, Acta Math. \textbf{155} (1985), no.~3-4,
  261--301.

\bibitem[dMG04]{Gonzalez2}
Mar\'ia del Mar~Gonz\'alez, \emph{Removability of singularities for a class of
  fully non-linear elliptic equations}, preprint, 2004.

\bibitem[dMG05]{Gonzalez1}
\bysame, \emph{Singular sets of a class of locally conformally flat manifolds},
  Duke Math. J. \textbf{129} (2005), no.~3, 551--572.

\bibitem[Eis97]{Eisenhart}
Luther~Pfahler Eisenhart, \emph{Riemannian geometry}, Princeton Landmarks in
  Mathematics, Princeton University Press, Princeton, NJ, 1997, Eighth
  printing, Princeton Paperbacks.

\bibitem[Eva82]{Evans}
Lawrence~C. Evans, \emph{Classical solutions of fully nonlinear, convex,
  second-order elliptic equations}, Comm. Pure Appl. Math. \textbf{35} (1982),
  no.~3, 333--363.

\bibitem[G{\.a}r59]{Garding}
Lars G{\.a}rding, \emph{An inequality for hyperbolic polynomials}, J. Math.
  Mech. \textbf{8} (1959), 957--965.

\bibitem[GT83]{GT}
David Gilbarg and Neil~S. Trudinger, \emph{Elliptic partial differential
  equations of second order}, second ed., Springer-Verlag, Berlin, 1983.

\bibitem[Gur93]{Gurskythesis}
Matthew~J. Gursky, \emph{Compactness of conformal metrics with integral bounds
  on curvature}, Duke Math. J. \textbf{72} (1993), no.~2, 339--367.

\bibitem[GV03]{GVJDG}
Matthew~J. Gursky and Jeff~A. Viaclovsky, \emph{A fully nonlinear equation on
  four-manifolds with positive scalar curvature}, J. Differential Geom.
  \textbf{63} (2003), no.~1, 131--154.

\bibitem[GV04]{GVAM}
Matthew Gursky and Jeff Viaclovsky, \emph{Volume comparison and the
  {$\sigma_k$}-{Y}amabe problem}, Adv. Math. \textbf{187} (2004), no.~2,
  447--487.

\bibitem[GVW03]{GVW}
Pengfei Guan, Jeff Viaclovsky, and Guofang Wang, \emph{Some properties of the
  {S}chouten tensor and applications to conformal geometry}, Trans. Amer. Math.
  Soc. \textbf{355} (2003), no.~3, 925--933 (electronic).

\bibitem[GW03a]{GuanWang2}
Pengfei Guan and Guofang Wang, \emph{A fully nonlinear conformal flow on
  locally conformally flat manifolds}, J. Reine Angew. Math. \textbf{557}
  (2003), 219--238.

\bibitem[GW03b]{GuanWang1}
\bysame, \emph{Local estimates for a class of fully nonlinear equations arising
  from conformal geometry}, Int. Math. Res. Not. (2003), no.~26, 1413--1432.

\bibitem[GW04a]{GuanWang3}
\bysame, \emph{Geometric inequalities on locally conformally flat manifolds},
  Duke Math. J. \textbf{124} (2004), 177--212.

\bibitem[GW04b]{GuanWangquot}
\bysame, \emph{Local gradient estimates for conformal quotient equations},
  preprint, 2004.

\bibitem[Har82]{Hartman}
Philip Hartman, \emph{Ordinary differential equations}, second ed.,
  Birkh\"auser Boston, Mass., 1982.

\bibitem[Kry83]{Krylov}
N.~V. Krylov, \emph{Boundedly inhomogeneous elliptic and parabolic equations in
  a domain}, Izv. Akad. Nauk SSSR Ser. Mat. \textbf{47} (1983), no.~1, 75--108.

\bibitem[Kui49]{Kuiper}
N.~H. Kuiper, \emph{On conformally-flat spaces in the large}, Ann. of Math. (2)
  \textbf{50} (1949), 916--924.

\bibitem[KW74]{KazdanWarner}
Jerry~L. Kazdan and F.~W. Warner, \emph{Curvature functions for compact
  {$2$}-manifolds}, Ann. of Math. (2) \textbf{99} (1974), 14--47.

\bibitem[Li89]{Yanyan2}
Yan~Yan Li, \emph{Degree theory for second order nonlinear elliptic operators
  and its applications}, Comm. Partial Differential Equations \textbf{14}
  (1989), no.~11, 1541--1578.

\bibitem[Lic58]{Lichnerowicz}
Andr{\'e} Lichnerowicz, \emph{G\'eom\'etrie des groupes de transformations},
  Travaux et Recherches Math\'ematiques, III. Dunod, Paris, 1958.

\bibitem[LL03]{LiLi2}
Aobing Li and Yan~Yan Li, \emph{On some conformally invariant fully nonlinear
  equations}, Comm. Pure Appl. Math. \textbf{56} (2003), no.~10, 1416--1464.

\bibitem[LP87]{LeeandParker}
John~M. Lee and Thomas~H. Parker, \emph{The {Y}amabe problem}, Bull. Amer.
  Math. Soc. (N.S.) \textbf{17} (1987), no.~1, 37--91.

\bibitem[MM03]{Mantegazza}
Carlo Mantegazza and Andrea~Carlo Mennucci, \emph{Hamilton-{J}acobi equations
  and distance functions on {R}iemannian manifolds}, Appl. Math. Optim.
  \textbf{47} (2003), no.~1, 1--25.

\bibitem[Pet98]{Petersen}
Peter Petersen, \emph{Riemannian geometry}, Graduate Texts in Mathematics, vol.
  171, Springer-Verlag, New York, 1998.

\bibitem[Sch89]{Schoen1}
Richard~M. Schoen, \emph{Variational theory for the total scalar curvature
  functional for {R}iemannian metrics and related topics}, Topics in calculus
  of variations (Montecatini Terme, 1987), Lecture Notes in Math., vol. 1365,
  Springer, Berlin, 1989, pp.~120--154.

\bibitem[Via00a]{Jeff1}
Jeff~A. Viaclovsky, \emph{Conformal geometry, contact geometry, and the
  calculus of variations}, Duke Math. J. \textbf{101} (2000), no.~2, 283--316.

\bibitem[Via00b]{Jeff3}
\bysame, \emph{Some fully nonlinear equations in conformal geometry},
  Differential equations and mathematical physics (Birmingham, AL, 1999), Amer.
  Math. Soc., Providence, RI, 2000, pp.~425--433.

\bibitem[Via02]{Jeff2}
\bysame, \emph{Estimates and existence results for some fully nonlinear
  elliptic equations on {R}iemannian manifolds}, Comm. Anal. Geom. \textbf{10}
  (2002), no.~4, 815--846.

\end{thebibliography}
\end{document}